\documentclass{article}
\usepackage{graphicx}
\usepackage{amssymb,amsmath}
\def\proof{\noindent\hspace{2em}{\itshape Proof: }}
\def\QEDclosed{\mbox{\rule[0pt]{1.3ex}{1.3ex}}} 

\def\QED{\QEDclosed}
\def\endproof{\hspace*{\fill}~\QED\par\endtrivlist\unskip}

\usepackage{color}
\usepackage{comment,xspace}
\usepackage{cite}
\newtheorem{theorem}{Theorem}
\newtheorem{definition}[theorem]{Definition}
\newtheorem{proposition}[theorem]{Proposition}
\newtheorem{lemma}[theorem]{Lemma}

\newtheorem{remark}[theorem]{Remark}

\def\reals{\mathbb{R}}
\def\naturals{\mathbb{N}}
\def\card{\mathrm{card}}
\def\vers{\mathrm{vers}}
\def\diam{\mathrm{diam}}

\begin{document}
\title{Efficient routing of multiple vehicles with no communications}
\author{Alessandro Arsie and Emilio Frazzoli \footnote{Research supported by National Science Foundation, and Air Force Office of Scientific Research, CCR-0325716, and F49620-02-1-0325, respectively.}
\\ Laboratory for Information and Decision Systems,\\
Massachusetts Institute of Technology, Cambridge, MA\\
arsie@ucla.edu, frazzoli@mit.edu}


\maketitle
\thispagestyle{empty}
\pagestyle{empty}

\begin{abstract}
In this paper we consider a class of dynamic vehicle routing problems, in which a number of mobile agents in the plane must visit target points generated over time by a stochastic process. It is desired to design motion coordination strategies in order to minimize the expected time between the appearance of a target point and the time it is visited by one of the agents. We propose control strategies that, while making minimal or no assumptions on communications between agents, provide the same level of steady-state performance achieved by the best known decentralized strategies. In other words, we demonstrate that inter-agent communication does not improve the efficiency of
such systems, but merely affects the rate of convergence to the steady state. Furthermore, the proposed strategies do not rely on the knowledge of the details of the underlying stochastic process. Finally, we show that our proposed strategies provide an efficient, pure Nash equilibrium in a game theoretic formulation of the problem, in which each agent's objective is to maximize the number of targets it visits. Simulation results are presented and discussed.
\end{abstract}

\section{Introduction}
A  very active research area today addresses coordination of several mobile agents: groups of autonomous robots and large-scale mobile networks are being considered for a broad class of applications, ranging from environmental monitoring, to search and rescue operations, and national security.

An area of particular interest is concerned with the generation of efficient cooperative strategies for several mobile agents to move through a certain number of given target points, possibly avoiding obstacles or threats~\cite{Beard.McLain.ea:02,Richards.Bellingham.ea:02,Schumacher.Chandler.ea:03,Earl.DAndrea:05,Li.Cassandras:06}.  Trajectory efficiency in these cases is understood in terms of cost for the agents: in other words, efficient trajectories minimize the total path length, the time needed to complete the task, or the fuel/energy expenditure.  A related problem has been investigated as the Weapon-Target Assignment (WTA) problem, in which mobile agents are allowed to team up in order to enhance the probability of a favorable outcome in a target engagement~\cite{Murphey:99,Arslan.Shamma:06}. In this setup, targets locations are known and an assignment strategy is sought that maximizes the global success rate.
In a biological setting, the closest parallel to many of these problems is the development of foraging  strategies, and of territorial vs. gregarious behaviors~\cite{Tanemura.Hasegawa:80}, in which individuals choose to identify and possibly defend a hunting ground.

In this paper we consider a class of cooperative motion coordination
problems, to which we can refer as {\em dynamic vehicle routing}, in which service requests are not known a priori, but are dynamically generated over time by a stochastic process in a geographic region of interest. Each service request is associated to a target point in the plane, and is fulfilled when one of a team of mobile agents visits that point. For example, service requests can be thought of as threats to be investigated in a surveillance application, events to be measured in an environmental monitoring scenario, and as information packets to be picked up and delivered to a user in a wireless  sensor network. It is desired to design a control strategy for the mobile agents that  provably minimizes the expected waiting time between the issuance of a service request and its fulfillment. In other words, our focus is on the quality of service as perceived by the ``end user," rather than, for example, fuel economies achieved by the mobile agents. Similar problems were also considered in~\cite{Psaraftis:88,Bertsimas.vanRyzin:91}, and decentralized strategies were presented in~\cite{Frazzoli.Bullo:CDC04}. This problem has connections to the Persistent Area Denial (PAD) and area coverage problems discussed, e.g., in~\cite{Schumacher.Chandler.ea:03,Liu.Cruz.ea:04,Cortes.Martinez.ea:04,Moore.Passino:05}.

A common theme in cooperative control is the investigation of the effects of different communication and information sharing protocols on the system performance. Clearly, the ability to access more information at each single agent can not decrease the performance level; hence, it is commonly believed that by providing better communication among agents will improve the system's performance.   In this paper, we prove that there are certain dynamic vehicle routing problems which can, in fact, be solved (almost) optimally without any explicit communication between agents; in other words, the no-communication constraint in such cases is not binding, and does not limit the steady-state performance. The main contribution of this paper is the introduction of a motion coordination strategy that does not require any explicit communication between agents, while achieving provably optimal performance in certain conditions.

The paper is structured as follows: in Section \ref{setup} we set up and formulate the problem we investigate in the paper.
In Section \ref{algorithm} we introduce the proposed solution algorithms, and discuss their characteristics. Section \ref{prova} is the technical core of the paper, in which we prove the convergence of the performance provided by the proposed algorithms to a critical point (either a local minimum or a saddle point) of the global performance function. Moreover, we show that any optimal configuration corresponds to a class of tessellations of the plane that we call Median Voronoi Tessellations.
Section \ref{games} is devoted to a game-theoretic interpretation of our result in which the agents are modeled as rational autonomous decision makers trying to maximize their own utility function. We prove that, following the policy prescribed by our algorithm, the agents reach an efficient pure Nash equilibrium which can be anyway suboptimal with respect to the global utility function (in this case the expected time of service).
In Section \ref{simulations} we present some numerical results, while Section \ref{conclusions} is dedicated to final remarks and further extensions of this line of research.

\section{Problem Formulation}\label{setup}
Let $\Omega \subset \reals^2$ be a convex
domain on the plane, with non-empty interior; we will refer to $\Omega$ as the {\em workspace}.  A stochastic process generates {\em service requests} over time, which are associated to points in $\Omega$; these points are also called {\em targets}. The process generating service requests is modeled as a spatio-temporal Poisson point process, with temporal intensity $\lambda >0$, and an absolutely continuous  spatial distribution described by the density function $\varphi: \Omega \rightarrow \reals_+$, with bounded and convex support within $\Omega$ (i.e., $\varphi(q)>0 \Leftrightarrow q \in \mathcal{Q} \subseteq \Omega$, with $\mathcal{Q}$ bounded and convex). The spatial density function $\varphi$ is normalized in such a way that $\int_\Omega \varphi(q) \; dq  = 1$. Both $\lambda$ and $\varphi$ are not necessarily known.

A spatio-temporal Poisson point process is a collection of functions $\{\mathcal{P}: \overline\reals_+ \rightarrow 2^\Omega\}$ such that, for any $t > 0$,  $\mathcal{P}(t)$ is a random collection of points in $\Omega$, representing the service requests generated in the time interval $[0, t)$, and such that
\begin{itemize}
\item The total numbers of events generated in two disjoint time-space regions are {\em independent} random variables;
\item The total number of events occurring in an interval $[s,s+t)$  in a measurable set $\mathcal{S} \subseteq \Omega$
satisfies
$$ \mathrm{ Pr}\left[\mathrm{card}\left( (\mathcal{P}(s+t) - \mathcal{P}(s))
\cap \mathcal{S} \right) =k \right]=\frac{ \exp(-\lambda t \cdot \varphi(\mathcal{S}))(\lambda  t\cdot \varphi(\mathcal{S}))^k}{k!},$$
 where this must holds for any $k$ in $\naturals$ and
where $\varphi(\mathcal{S})$ is a shorthand for $\int_\mathcal{S} \varphi(q) \; dq$.
\end{itemize}
Each particular function $\mathcal{P}$ is a realization, or trajectory, of the Poisson point process.
A consequence of the properties defining Poisson processes is that the {\em expected} number of targets generated in a measurable region $S \subseteq \Omega$ during a time interval of length $\Delta t$ is given by: $$ \mathrm{E}[\mathrm{card}\left( (\mathcal{P}(t+\Delta t) - \mathcal{P}(t)) \cap \mathcal{S}\right)]=
\lambda \Delta t \cdot \varphi(\mathcal{S}).$$
Without loss of generality, we will identify service requests with targets points, and label them in order of generation; in other words, given two targets $e_i, e_j \in \mathcal{P}(t)$, with $i<j$, the service request associated with these target have been issued at times $t_i \le t_j \le t$ (since events are almost never generated concurrently, the inequalities are in fact strict almost surely).

A service request is fulfilled when one of $m$ mobile agents, modeled as point masses, moves to the target point associated with it; $m$ is a possibly large, but finite number. Let
$p(t) = \left(p_1(t), p_2 (t), \ldots, p_m(t)\right) \in \Omega^m$ be a vector describing the positions of the agents at time $t$. (We will tacitly use a similar notation throughout the paper). The agents are free to move, with bounded speed, within the  workspace $\Omega$; without loss of generality, we will assume that the maximum speed is unitary.
In other words,  the dynamics of the agents are described by differential equations of the form
\begin{equation}
\label{eq:model}
\frac{d\; p_i(t)}{dt} = u_i(t), \quad \mbox{ with } \|u_i(t)\| \le 1,\quad \forall t \ge 0, i \in \{1, \ldots, m\}.
\end{equation}
The agents are identical, and have unlimited range and target-servicing capability.

Let $\mathcal{B}_i(t) \subset \Omega$ indicate the set of targets serviced by the $i$-th  agent up to time $t$. (By convention, $\mathcal{B}_i(0) = \emptyset$, $i = 1, \ldots, m$). We will assume that $\mathcal{B}_i \cap \mathcal{B}_j = \emptyset$ if $i \neq j$, i.e., that service requests are fulfilled by at most one agent. (In the unlikely event that two or more agents visit a target at the same time, the target is arbitrarily assigned to one of them).

Let $\mathcal{D}: t \rightarrow 2^{\Omega}$ indicate (a realization of) the stochastic process obtained combining the service request generation process $\mathcal{P}$ and the  removal process caused by the agents servicing outstanding requests; in other words,
$$\mathcal{P}(t) = \mathcal{D}(t) \cup \mathcal{B}_1(t) \cup \ldots \cup \mathcal{B}_m(t), \qquad \mathcal{D}(t) \cap \mathcal{B}_i(t) = \emptyset, \; \forall i \in \{1, \ldots, m\}.$$
The random set $\mathcal{D}(t) \subset \Omega$ represents the {\em demand}, i.e., the service requests outstanding at time $t$; let $n(t) = \mathrm{card}(\mathcal{D}(t))$.

Our objective in this paper will be the design of motion coordination strategies that allow the mobile agents to fulfill service requests efficiently (we will make this more precise in the following). In particular, in this paper we will concentrate on  motion coordination strategies of the following two forms:
\begin{equation}
\label{eq:nocomm}
\pi_i: (p_i , \mathcal{B}_i, \mathcal{D}) \mapsto u_i,\qquad i \in \{1, \ldots, m\},
\end{equation}
and
\begin{equation}
\label{eq:sensorbased}
\pi_i: (p_1, \ldots, p_m , \mathcal{B}_i, \mathcal{D}) \mapsto u_i, \qquad i \in \{1, \ldots, m\}.
\end{equation}
An agent executing a control policy of the form \eqref{eq:nocomm} relies on the knowledge of its own current position, on a record of targets it has previously visited, and on the current demand.
In other words, such control policies do not need any explicit information exchange between agents; as such, we will refer to them as {\em no communication} ($\mathrm{nc}$) policies. Such policies are trivially decentralized.

On the other hand, an agent executing a control policy of the form \eqref{eq:sensorbased} can sense the current position of other agents, but still has information only on the targets itself visited in the past (i.e., does not know what, if any, targets have been visited by other agents).  We call these {\em sensor-based} ($\mathrm{sb}$) policies, to signify the fact that only factual information is exchanged between agents---as opposed to information related to intent and past history.  Note that both families of coordination policies rely, in principle, on the knowledge of the locations of all outstanding targets. (However, as we will see in the following, only local target sensing will be necessary in practice).

A policy $\pi = (\pi_1, \pi_2, \ldots, \pi_m)$ is said to be {\em stabilizing} if, under its effect, the expected number of outstanding targets does not diverge over time, i.e., if
\begin{equation}
\overline n_\pi = \lim_{t\rightarrow \infty} \mathrm{E}[ n(t) \| \dot p_i(t) = \pi_i(
p(t), \mathcal{B}_i(t), \mathcal{D}(t)), i \in \{1, \ldots, m\}] < \infty.
\end{equation}
Intuitively, a policy is stabilizing if the mobile agents are able to visit targets at a rate that is---on average---at least as fast as the rate at which new service requests are generated.

Let $T_j$ be the time elapsed between the issuance of the $j$-th service request, and the time it is fulfilled. If the system is stable,
then the following balance equation (also known as Little's formula~\cite{Little:61}) holds:
\begin{equation}
\label{eq:little}
\overline n_\pi  = \lambda \overline T_\pi,
\end{equation}
where $\overline{T}_\pi:=\lim_{j\rightarrow \infty} \mathrm{E}[T_j]$ is the system time under policy $\pi$, i.e., the expected time  a service request must wait before being fulfilled, given that  the mobile agents follow the strategy defined by $\pi$.  Note that  the system time $\overline T_\pi$ can be thought of as a measure of the quality of service, as perceived by the ``user" issuing the service requests.

At this point we can finally state our problem: we wish to devise a policy that is (i) stabilizing, and (ii) yields a quality of service (i.e., system time)  achieving, or approximating, the theoretical optimal performance given by
\begin{equation}
\label{eq:opt}
\overline T_\mathrm{opt} = \inf_{\pi \text{ stabilizing}} \overline T_\pi
\end{equation}

Centralized and decentralized strategies are known that optimize or approximate
\eqref{eq:opt} in a variety of cases of interest~\cite{Bertsimas.vanRyzin:91,Bertsimas.vanRyzin:93,Bertsimas.vanRyzin:93b,Frazzoli.Bullo:CDC04}. However, all such strategies rely either on a central authority with the ability to communicate to all agents, or on the exchange of certain information about each agent's strategy with other neighboring agents. In addition, these policies require the knowledge of the spatial distribution $\varphi$; decentralized versions of these implement versions of Lloyd's algorithm for vector quantization~\cite{Lloyd:82}.

In the remainder of this paper, we will investigate how the additional constraints posed on the exchange of information between agents by the models \eqref{eq:nocomm} and \eqref{eq:sensorbased} impact the achievable performance and quality of service. Remarkably, the  policies we will present do not rely on the knowledge of the spatial distribution $\varphi$, and are a generalized version of MacQueen's clustering algorithm~\cite{MacQueen:67}.

\section{Control policy description}
\label{algorithm}
\begin{figure}[tbh]
\centerline{\includegraphics[width=0.8\textwidth]{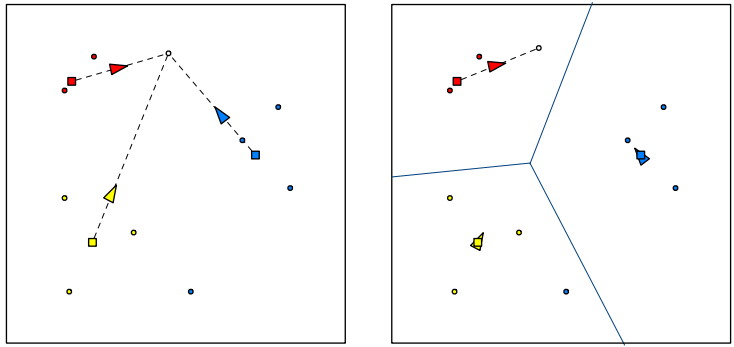}}
\caption{An illustration of the two control policies proposed in
Section \ref{algorithm}. While no targets are outstanding,
vehicles wait at the point that minimizes the average distance to
targets they have visited in the past; such points are depicted as
squares, while targets are circles and vehicles triangles. In the
no-communication policy, at the appearance of a new target, all
vehicles pursue it (left). In the sensor-based policy, only the
vehicle that is closest to the target will pursue it (right).}
\label{fig:algorithms}
\end{figure}

In this section, we introduce two control policies of the forms, respectively, \eqref{eq:nocomm} and \eqref{eq:sensorbased}. An illustration of the two policies is given in Figure \ref{fig:algorithms}.

\subsection{A control policy requiring no explicit communication}
Let us begin with an informal description of a policy $\pi_\mathrm{nc}$ requiring no explicit information exchange between agents. At any given time $t$, each agent computes its own control input according to the following rule:
\begin{enumerate}
\item If $\mathcal{D}(t)$ is not empty, move towards the nearest outstanding target.
\item If $\mathcal{D}(t)$ is empty, move towards the point minimizing the average distance to targets {\em serviced in the past} by each agent. If there is no unique minimizer, then move to the nearest one.
\end{enumerate}
In other words, we set
\begin{equation}
\pi_\mathrm{nc} (p_i(t), \mathcal{B}_i(t), \mathcal{D}(t)) =
\vers( F_\mathrm{nc}(p_i(t), \mathcal{B}_i(t), \mathcal{D}(t)) - p_i(t)),
\end{equation}
where
\begin{equation}
\label{eq:Fnc}
F_\mathrm{nc}(p_i,\mathcal{B}_i, \mathcal{D}) =
\left\{
\begin{array}{ll} \displaystyle \arg \min_{q \in \mathcal{D}} \|p_i - q\|, & \mbox{ if } \mathcal{D} \neq \emptyset, \\[10pt]
\displaystyle \arg \min_{q \in \Omega} \sum_{e \in \mathcal{B}_i} \|e - q\|, & \mbox{ otherwise, }
\end{array}
\right.
\end{equation}
$\|\cdot\|$ is the Euclidean norm, and
$$\vers(v) = \left\{ \begin{array}{ll} v/\|v\|, & \mbox{ if } v\neq 0,\\
0 & \mbox{ otherwise.}\end{array}\right.$$
The convex function $W: q \mapsto \sum_{e \in \mathcal{B}} \|q - e\| $, often called the (discrete) Weber function in the facility location literature~\cite{Drezner:95,Agarwal.Sharir:98} (modulo normalization by $\mathrm{card}(\mathcal{B})$), is not strictly convex only when the point set $\mathcal{B}$ is empty---in which case we set $W(\cdot)=0$ by convention--- or contains an even number of collinear points. In such cases, the minimizer nearest to $p_i$ in \eqref{eq:Fnc} is chosen. We will call the point $p^*_i(t) = F_\mathrm{nc}(\cdot, \mathcal{B}_i(t), \emptyset)$ the {\em reference point} for the $i$-th agent at time $t$.

In the $\pi_\mathrm{nc}$ policy, whenever one or more service requests are outstanding, all agents will be pursuing a target; in particular, when only one service request is outstanding, all agents will move towards it. When the demand queue is empty, agents will either (i) stop at the current location, if they  have visited no targets yet, or (ii) move to their  reference point, as determined by the set of targets previously visited.


\subsection{A sensor-based control policy}
\label{sensorbc}
The control strategy in the previous section can be modified to include
information on the current position of other agents, if available (e.g., through on-board sensors).  In order to present the new policy, indicate with $\mathcal{V}(p) = \{\mathcal{V}_1(p), \mathcal{V}_2(p), \ldots, \mathcal{V}_m(p)\}$ the Voronoi partition of the workspace $\Omega$, defined as:
\begin{equation}
\label{eq:Voronoi}
\mathcal{V}_i(p) = \{ q \in \Omega: \|q-p_i\| \le \|q-p_j\|, \forall j = 1\ldots m \}.
\end{equation}

As long as an agent has never visited any target, i.e., as long as $\mathcal{B}_i(t) = \emptyset$, it executes the $\pi_\mathrm{nc}$ policy. Once an agent has visited at least one target, it computes its own control input according to the following rule:
\begin{enumerate}
\item If $\mathcal{D}(t) \cap \mathcal{V}_i(t)$ are not empty, move towards the nearest outstanding target in the agent's own Voronoi region.
\item If $\mathcal{D}(t) \cap \mathcal{V}_i(t)$ is empty, move towards the point minimizing the average distance to targets in $\mathcal{B}_i(t)$.  If there is no unique minimizer, then move to the nearest one.
\end{enumerate}
In other words, we set
\begin{equation}
\pi_\mathrm{sb} (p(t), \mathcal{B}_i(t), \mathcal{D}(t)) =
\vers( F_\mathrm{sb}(p(t), \mathcal{B}_i(t), \mathcal{D}(t)) - p_i(t)),
\end{equation}
where
\begin{equation}
\label{eq:pi_nc}
F_\mathrm{sb}(p,\mathcal{B}_i, \mathcal{D}) =
\left\{
\begin{array}{ll}
\displaystyle \arg \min_{q \in \mathcal{D}} \|p_i - q\|, & \mbox{ if } \mathcal{D} \cap \mathcal{V}_i \neq \emptyset, \mbox{ and } \mathcal{B}_i = \emptyset\\[10pt]
\displaystyle \arg \min_{q \in \mathcal{D} \cap \mathcal{V}_i(p)} \|p_i - q\|, & \mbox{ if } \mathcal{D} \cap \mathcal{V}_i \neq \emptyset, \mbox{ and } \mathcal{B}_i \neq \emptyset\\[10pt]
\displaystyle \arg \min_{q \in \Omega} \sum_{e \in \mathcal{B}_i} \|e - q\|, & \mbox{ otherwise. }
\end{array}
\right.
\end{equation}

In the $\pi_\mathrm{sb}$ policy, at most one agent will be pursuing a given target, at any time after an initial transient that terminates when all agents have visited at least one target each. The agents' behavior when no oustanding targets are available in their Voronoi region is similar to that determined by the $\pi_\mathrm{nc}$ policy previously discussed, i.e., they move to their reference point, determined by previously visited targets.

\begin{figure}[tbh]
\centerline{\includegraphics[width=0.3\columnwidth]{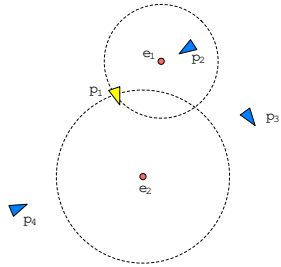}}
\caption{Implicit computation of Voronoi regions: Even though
target $e_1$ is the nearest target to $p_1$, it is not in the
Voronoi region of the 1st agent. In fact, the circle of radius
$\|e_1-p_1\|$ centered at $e_1$ contains $p_2$, and the 2nd agent
is closer to $e_1$. However, the circle of radius $\|e_2 - p_1$
centered at $e_2$ does not contain any other agent, ensuring that
$e_2$ is in the Voronoi region generated by $p_1$.}
\label{fig:radius}
\end{figure}

\begin{remark} While we introduced Voronoi partitions in the definition of the control policy, the explicit computation of each agent's Voronoi region is not necessary.  In fact, each agent only needs to check whether it is the closest agent to a given target or not. In order to check whether a target point $q$ is in the Voronoi region of the $i$-th agent, it is necessary to know the current position only of agents within a circle or radius $\|p_i-q\|$ centered at $q$ (see Figure \ref{fig:radius}). For example, if such circle is empty, then $q$ is certainly in $\mathcal{V}_i$; if the circle is not empty, distances of the agents within it to the target must be compared. This provides a degree of spatial decentralization---with respect to other agents---that is even stronger than that provided by restricting communications to agents sharing a boundary in a Voronoi partition (i.e., neighboring agents in the Delaunay graph, dual to the partition \eqref{eq:Voronoi}).
\end{remark}

\begin{remark} The sensor-based policy is more efficient than the no-communication policy in terms of the length of the path traveled by each agent, since there is no duplication of effort as several agents pursue the same target. However, in terms of ``quality of service," we will show that there is no difference between the two policies, for low target generation rates. Numerical results show that the sensor-based policy is more efficient in a broader range of target generation rates, and in fact provides almost optimal performance both in light and heavy load conditions.
\end{remark}

\section{Performance analysis in light load}\label{prova}
In this section we analyze the performance  of the control policies proposed
in the previous section. In particular, we concentrate our investigation
on the light load case, in which the target generation rate is very small, i.e., as $\lambda \rightarrow 0^+$. This will allow us to prove analytically certain interesting and perhaps surprising characteristics of the proposed policies. The performance analysis in the general case is more difficult; we will discuss the results of numerical investigation in Section \ref{simulations}, but no analytical results are available at this time.

\subsection{Overview of the system behavior in the light load regime}
Before starting a formal analysis, let us summarize the key characteristics of the agents' behavior in light load, i.e., for small values of  $\lambda$.
\begin{enumerate}
\item At the initial time the $m$ agents are assumed to be deployed in general position in $\Omega$, and the demand queue is empty, $\mathcal{D}(0) = \emptyset$.
\item The agents do not move until the first service request appears.
At that time, if the policy $\pi_\mathrm{nc}$ is used, all agents will start moving towards the first target. If the sensor-based policy $\pi_\mathrm{sb}$ is used, only the closest agent will move towards the target.
\item As soon as one agent reaches the target, all agents start moving towards their current reference point, and the process continues.
\end{enumerate}

For small $\lambda$, with high probability (i.e., with probability approaching 1 as $\lambda \to 0$) at most one service request is outstanding at any given time. In other words, new service requests are generated so rarely that most of the time agents will be able to reach a target and return to their reference point before a new service request is issued.


Consider the $j$-th service request, generated at time $t_j$. Assuming that at $t_j$ all agents are at their reference position, the expected system time $T_j$ can be computed as
$$T_j = \int_\Omega \min_{i=1, \ldots, m} \|p^*_i(t_j) - q\| \; \varphi(q) dq.$$
Assume for now that the sequences $\{p^*_i(t_j): j \in \naturals\}$ converge,
and let  $$\lim_{j \rightarrow \infty} p^*_i(t_j) = \hat{p}^*_i.$$
Note that $\hat p^*_i$ is a random variable, the value of which depends in general on the particular realization of the target generation process.
If all service requests are generated with the agents at their reference position, the average service time (for small $\lambda$) can be evaluated as
\begin{equation}
\overline T_\mathrm{nc} = \overline T_\mathrm{sb} =  \int_\Omega \min_{i=1, \ldots, m} \left\|  \hat p^*_i - q\right\| \; \varphi(q) dq  = \sum_{i=1}^m \int_{\mathcal{V}_i(p^*)} \|\hat p^*_i - q\| \; \varphi(q) dq.
\end{equation}
Since the system time depends on the random variable $\hat p^*= (\hat p_1^*, \dots,\hat p_m^*)$, it is itself a random variable. The function appearing on the right hand side of the above equation, relating the system time to the asymptotic location of reference points, is called the continuous multi-median function~\cite{Drezner:95}. This function admits a global minimum (in general not unique) for all non-singular density functions $\varphi$, and in fact it is known~\cite{Bertsimas.vanRyzin:91} that the optimal performance in terms of system time is given by
\begin{equation}
\label{eq:Topt}
\overline T_\mathrm{opt} = \min_{p \in {\Omega}^m} \sum_{i=1}^m \int_{\mathcal{V}_i(p)} \|p_i - q\| \; \varphi(q) dq.
\end{equation}

In the following, we will investigate the convergence of the reference points as new targets are generated, in order to draw conclusions about the average system time $\overline T$ in light load. In particular, we will prove not only that the reference points converge with high probability (as $\lambda \to 0$) to a local critical point (more precisely, either local minima or saddle points) for the average system time, but also that the limiting reference points $\hat p^*$ are {\em generalized medians} of their respective Voronoi regions, where
 \begin{definition}[Generalized median]
The  {\em generalized median} of a set $\mathcal{S}\subset \mathbb{R}^n$ with respect to a density function $\varphi:\mathcal{S} \to \overline{\mathbb{R}}_+$ is defined as
$$\overline{p}:=\arg\min_{p \in \mathbb{R}^n} \int_{\mathcal{S}}\|p-q\|\varphi(q)\; dq.$$
\end{definition}
We call the resulting Voronoi tessellation {\em Median Voronoi Tessellation} (MVT for short), in analogy with what is done with Centroidal Voronoi Tessellations.
A formal definition is as follows:
\begin{definition}[Median Voronoi Tessellation]\label{MVT}
A Voronoi tessellation $\mathcal{V}(p)= \{\mathcal{V}_1(p), \ldots, \mathcal{V}_m(p)\}$ of a set $\mathcal{S} \subset \mathbb{R}^n$ is said a {\em Median Voronoi Tessellation} of $\mathcal{S}$ with respect to the density function $\varphi$  if the ordered set of generators $p$ is equal to the ordered set of generalized medians of the sets in $\mathcal{V}(p)$ with respect to $\varphi$, i.e., if
$$p_i = \arg \min_{s \in \mathbb{R}^n} \int_{\mathcal{V}_i(p)} \|s-q\| \varphi(q)\; dq, \qquad \forall i \in \{1, \ldots, m\}.$$
\end{definition}

Since the proof builds on a number of intermediate results, we provide an outline of the argument as a convenience to the reader.
\begin{enumerate}
\item First we prove that the reference point of any agents that visits an unbounded number of targets over time converges almost surely.
\item Second, we prove that, if $m \ge 1$ agents visit an unbounded number of targets over time, their reference points will converge to the generators of a MVT almost surely, as long as agents are able to return to their reference point infinitely often.
\item Third, we prove that all agents will visit an unbounded number of targets (this corresponds to a property of distributed algorithms that is often called {\em fairness} in computer science).
\item Finally, we prove that agents are  able to return to their reference point infinitely often with high probability as $\lambda \to 0^+$.
\end{enumerate}
Combining these steps, together with \eqref{eq:opt}, will allow us to state that the reference points converge to a local critical point of the system  time, with high probability as $\lambda \to 0^+$.

\subsection{Convergence of reference points}
Let us consider an agent $i$, such that
$$\lim_{t\rightarrow \infty} \card(\mathcal{B}_i(t)) = \infty,$$
i.e., an agent that services an unbounded number of requests over time.
Since the number of agents $m$ is finite, and the expected number of targets generated over a time interval $[0, t)$ is proportional to $t$, at least one such agent will always exist. In the remainder of this section, we will drop the subscript $i$, since we will consider only this agent, effectively ignoring all others for the time being.

For any finite $t$, the set $\mathcal{B}(t)$ will contain a finite number of points. Assuming that $\mathcal{B}(t)$ contains at least three non-collinear points, the discrete Weber function $p \mapsto \sum_{q \in \mathcal{B}(t)} \|p - q\|$ is strictly convex, and has a unique optimizer
$p^*(t) = \arg\min_{p \in \Omega} \sum_{q \in \mathcal{B}(t)} \|p - q\|$.
The optimal point $p^*(t)$ is called the Fermat-Torricelli (FT) point---or the Weber point in the location optimization literature---associated with the set  $\mathcal{B}(t)$; see~\cite{Chandrasekaran.Tamir:90,Wesolowsky:93,Agarwal.Sharir:98} for a historical review of the problem and for solution algorithms.

It is known that the FT point is unique and algebraic for any set of non-collinear points. While there are no general analytic solutions for the location of the FT point associated to more than 4 points, numerical solutions can be easily constructed relying on the convexity of the Weber function, and on the fact that is is differentiable for all points not in $\mathcal{B}$. Polynomial-time approximation algorithms are also available (see, e.g., \cite{Fekete.Mitchell.ea:00,Carmi.Har-Peled.Katz:05}). Remarkably, a simple mechanical device can be constructed to solve the problem, based on the so-called Varignon frame, as follows. Holes are drilled on a horizontal board, at locations corresponding to the points in $\mathcal{B}$.
A string attached to a unit mass is passed through each of these holes, and all strings are tied together at one end. The point reached by the knot at equilibrium is a FT point for $\mathcal{B}$.

Some useful properties of FT points are summarized below.
If there is a $q_0 \in \mathcal{B}$ is such that
\begin{equation}
\label{eq:absorbing}
\left\| \sum_{q \in \mathcal{B}\setminus q_0} \vers(q_0-q) \right\| \le 1
\end{equation}
then $p^* = q_0$ is a FT point for $\mathcal{B}$. If no point in $\mathcal{B}$ satisfies such condition, then the FT point $p^*$ can be found as a solution of the
following equation:
\begin{equation}
\label{eq:unitvectors}
\sum_{q \in \mathcal{B}} \vers(p^*-q) = 0.
\end{equation}
In other words, $p^*$ is simply the point in the plane at which the sum of unit vectors starting from it and directed to each of the points in $\mathcal{B}$ is equal to zero; this point is unique if $\mathcal{B}$ contains non-collinear points. Clearly, the FT point is in the  convex hull of $\mathcal{B}$.

Note that the points in $\mathcal{B}(t)$ are randomly sampled from an unknown absolutely continuous distribution, described by a spatial density function $\tilde \varphi$---which is not necessarily the same as $\varphi$, and in general is time-varying, depending on the past actions of all agents in the system.  Even though $\tilde \varphi$ is not known, it can be expressed as
$$\varphi(q,t) = \left\{\begin{array}{ll}
\displaystyle \frac{\varphi(q)}{\int_{\mathcal{I}(t)} \varphi(q) \; dq} & \mbox{ if } q \in \mathcal{I}(t)\\
0 & \mbox{ otherwise,}\end{array}\right.$$ for some convex set $\mathcal{I}(t)$ containing $p(t)$. (In practical terms, such set will be the Voronoi region generated by $p(t)$).

The function $t \mapsto p^*(t)$ is piecewise constant, i.e., it changes value at the discrete time instants $\{t_j: j \in \naturals\}$ at which  the agent visits new targets. As a consequence, we can concentrate on the sequence $\{p^*(t_j): j \in \naturals\}$, and study its convergence.

\begin{definition}
For any $t>0$, let the {\em solution set} C(t) be defined as
$$C(t):=\left\{p\in \Omega : \left\| \sum_{q \in \mathcal{B}(t)} \vers(p-q) \right\| \leq 1\right \}.$$
\end{definition}
\begin{figure}[tbh]
\centerline{\includegraphics[width=0.6\textwidth]{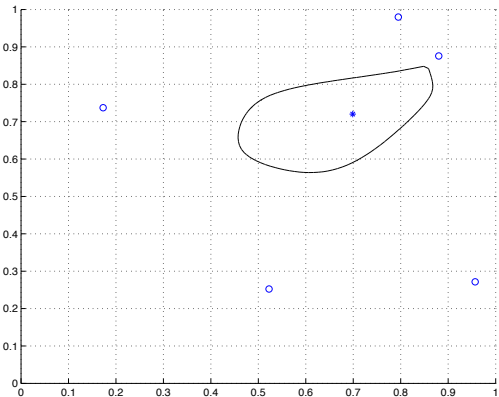}}
\caption{Example of a Fermat-Torricelli point (star) and solution
set corresponding to five target points (circles). Upon the
addition of an arbitrarily chosen sixth target point, the
Fermat-Torricelli is guaranteed to remain within the region
bounded by the curve.} \label{fig:solutionset}
\end{figure}
An example of such set is shown in Figure \ref{fig:solutionset}. The reason for introducing such solution sets  is that they have quite remarkable properties as shown by the following

\begin{proposition}
For any $j \in \naturals$, $p^*(t_{j+1}) \in C(t_j)$.
More specifically, if $e_{j+1} \notin C(t_j)$ (i.e., the target point associated to the $j$-th service request is outside the solution set) then the FT point $p^*(t_{j+1})$ is on the boundary of $C_i$.  If $e_{j+1}\in C(t_j)$, then $p^*(t_{j+1})=e_{j+1}$.
\end{proposition}
\proof If $e_{j+1}$ lies outside $C(t_j)$, we search for $p^*(t_{j+1})$ as the solution of the equation
$$\sum_{q \in \mathcal{B}(t_j)} \vers(p-q) + \vers(p-e_{j+1})  = 0,$$ from which it turns out immediately
 $$\left\|\sum_{q \in \mathcal{B}(t_j)} \vers(p-q)\right\| =\left\| -\vers(p-e_{j+1})\right\|=1,$$
thus $p^*(t_{j+1})\in \partial C(t_j)$. Notice that is is not true in general that the solution $p^*(t_{j+1})$ will lie on the line connecting $p^*(t_j)$ with the new target $e_{j+1}$.  In the other case, if $e_{j+1}$ lies in $C(t_j)$, then it satisfies condition \eqref{eq:absorbing}, and is the new FT point.
\endproof

Now, in order to prove that the $\{p^*(t_j)\}_{j\in\naturals}$ converges to a point $\hat{p}^*$, we will prove that the diameter of the solution set
$C(t_j)$ vanishes almost surely as $j$ tends to infinity.
%
%
%
%
First we need the following result.

\begin{figure}[htb]
\centerline{\includegraphics[width=0.5\textwidth]{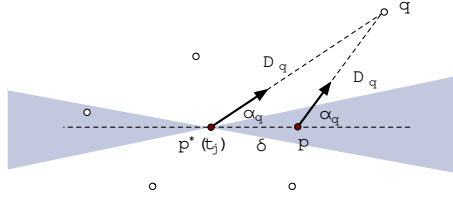}}
\caption{Geometric constructions in the proof of Proposition
\ref{prop1}.} \label{fig:proof_alpha}
\end{figure}

  \begin{proposition}\label{convex}
 If $\mathcal{Q} = \mathrm{Supp}(\varphi)$ is convex with non-empty interior, then $p^*(t) \in \mathrm{int} (\mathcal{Q})$ almost surely, for all $t$ such that $\mathcal{B}(t) \neq \emptyset$ .
 \end{proposition}
 \proof For any non-empty target set $\mathcal{B}(t)$, the FT point lies within the convex hull of $\mathcal{B}(t)$. All points in the set $\mathcal{B}(t)$ are contained within the interior of $\mathcal{Q}$ with probability one, since the boundary of $\mathcal{Q}$ is a set of measure zero. Since $\mathrm{int}(\mathcal{Q})$ is convex, and $\mathcal{B}(t) \subset \mathrm{int}(\mathcal{Q})$ almost surely,
$p^*(t) \in \mathrm{co}(\mathcal{B}(t)) \subset \mathrm{int}(\mathcal{Q}),$
almost surely.
\endproof
\begin{proposition}\label{prop1}
If the support of $\varphi$ is convex and bounded,
$$\lim_{j\rightarrow\infty}\mathrm{diam}(C(t_j))=0, \qquad\mathrm{a.s.}$$
\end{proposition}
\proof
Consider a generic point $p \in C(t_j)$, and let $\delta = p-p^*(t_j)$, $$\alpha_q = \arccos \left[\mathrm{vers}(p-p^*(t_j)) \cdot \mathrm{vers}(q-p^*(t_j))\right]$$,  $$\alpha'_q = \arccos \left[\mathrm{vers}(p-p^*(t_j)) \cdot \mathrm{vers}(q-p)\right],$$ see Figure \ref{fig:proof_alpha}.

Since $p \in C(t_j)$, the magnitude of the sum of unit vectors $\sum_{q \in \mathcal{B}(t_j)}\vers(p-q)$ is no more than one, and the following inequality is true:
\begin{multline}
\label{eq:cosalphaprime}
\left|\left(\sum_{q \in \mathcal{B}(t_j)} \vers(p-q)\right) \cdot \vers(p-p^*(t_j))\right| \\= \left|\sum_{q\in \mathcal{B}(t_j)} \left(\vers(p-q)\cdot \vers(p-p^*(t_j)) \right) \right| \\= \left| \sum_{q \in \mathcal{B}(t_j)}
\cos(\alpha'_q)\right|\le 1.
\end{multline}

Using elementary planar geometry, we obtain that
$$\alpha'_q - \alpha_q \ge \sin(\alpha'_q - \alpha_q) = \frac{\delta \sin (\alpha_q) }{\|q-p\|}.$$
Pick a small angle $0 < \alpha_\mathrm{min} < \pi/2$, and let
$$\mathcal{B}_{\alpha_\mathrm{min}}(t_j) = \{q \in \mathrm{B}(t_j): \sin (\alpha_q) \ge \sin(\alpha_\mathrm{min})\}.$$
(In other words, $\mathcal{B}_{\alpha_\mathrm{min}}(t)$ contains all points in $\mathcal{B}$ that are not in a conical region of half-width $\alpha_\mathrm{min}$, as shown in Figure \ref{fig:proof_alpha}).
For all $q \in \mathcal{B}(t_j)$, $\cos(\alpha'_q) \le \cos(\alpha_q)$; moreover, for all $q \in \mathcal{B}_{\alpha_\mathrm{min}}$,
$$\cos(\alpha'_q) \le \cos(\alpha_q) - \sin(\alpha_\mathrm{min}) (\alpha'_q - \alpha_q) \le$$
$$\le \cos(\alpha_q) -  \frac{\delta\sin(\alpha_\mathrm{min})^2}{\|q - p\|}
\le \cos(\alpha_q) -  \frac{\delta\sin(\alpha_\mathrm{min})^2}{\mathrm{diam}(\mathcal{Q})+\delta}.$$
Hence, summing over all $q \in \mathcal{B}(t_j)$, we get: 
\begin{equation}
\label{eq:cardB}
\sum_{q \in \mathcal{B}(t_j)} \cos (\alpha'_q) \le \sum_{q \in \mathcal{B}(t_j)} \cos (\alpha_q) - \sum_{q \in \mathcal{B}_{\alpha_\mathrm{min}}(t_j)}  \frac{\delta\sin(\alpha_\mathrm{min})^2}{\mathrm{diam}(\mathcal{Q})+\delta} 
\end{equation}
Observe now that in any case $$\left|  \sum_{q \in \mathcal{B}(t_j)} \cos (\alpha_q) \right| \leq 1,$$
(it is zero in case $p^*(t_j)\notin \mathcal{B}(t_j)$, and bounded in absolute value by one if
 $p^*(t_j)\in \mathcal{B}(t_j)$).
 Therefore, rearranging equation \eqref{eq:cardB}:
 $$\mathrm{card}(\mathcal{B}_{\alpha_\mathrm{min}}(t_j)) \frac{\delta\sin(\alpha_\mathrm{min})^2}{\mathrm{diam}(\mathcal{Q})+\delta}\leq \sum_{q \in \mathcal{B}(t_j)} \cos (\alpha_q) - \sum_{q \in \mathcal{B}(t_j)} \cos (\alpha'_q) \leq 2$$
Solving this inequality with respect to $\delta$ we get:
\begin{equation}
\label{eq:delta_bound}
\delta \le \frac{2\,\diam (\mathcal{Q})} {\mathrm{card}(\mathcal{B}_{\alpha_\mathrm{min}}(t_j)) \sin(\alpha_\mathrm{min})^2 -2}.
\end{equation}

Since (i) $\alpha_\mathrm{min}$ is a positive constant, (ii) $p^*(t_j)$ is in the interior of $\mathcal{Q}$, and (iii) $\lim_{j\to\infty}\card(\mathcal{B}(t_j)) = +\infty$, the right hand side of \eqref{eq:delta_bound} converges to zero with probability one. Since the bound holds for all points $p \in C(t_j)$, for all $j\in \mathbb{N}$, the claim follows.
\endproof

In the previous proposition, we have proven that $\|p^*(t_{j+1})-p^*(t_j)\|$ tends to zero a.s. as $j\rightarrow \infty$, under some natural assumptions on the distribution $\varphi$ and its support.
Unfortunately this is not sufficient to prove that  the sequence $\{p^*(t_j)\}_{t_j\in \naturals}$ is Cauchy; convergence of the sequence is however ensured by the following
\begin{proposition}\label{thm:convergence}
The sequence $\{p^*(t_j)\}_{j\in \naturals}$ converges almost surely.
\end{proposition}
\proof Since the sequence  $\{p^*(t_j)\}_{j\in \mathbb{N}}$ takes value in a compact set (the closure of $\mathrm{Supp}(\varphi)$), by the Bolzano-Weirstrass theorem there exists a subsequence converging to a limit point $\hat{p}^*$ in the compact set. Construct from  $\{p^*(t_j)\}_{j\in \mathbb{N}}$ a maximal subsequence converging to $\hat{p}^*$, and call $J$ the set of indices of this maximal subsequence. If the original sequence $\{p^*(t_j)\}_{j\in \mathbb{N}}$ is not converging to $\hat{p}^*$, then there exists an $L>0$ such that $\|p^*(t_j)-\hat{p}^*\|\geq L$, for any $j \in \naturals \setminus J$, and this set of indices is unbounded. Take $\epsilon:=L/3>0$. We have that $\|p^*(t_{k-1})-\hat{p}^*\|\leq \epsilon$ for any sufficiently large $(k-1)\in J $; moreover, $\|p^*(t_{k-1})-p^*(t_k)\|\leq \epsilon$, a.s. by Proposition (\ref{prop1}). Choose a sufficiently large $(k-1)\in J$, such that $k\in\mathbb{N}- J$ (this is always possible since the complementary set of $J$ is unbounded by the assumption of non-convergence). But now, we have
$$ L\leq \|p^*(t_k)-\hat{p}^*\|\leq \|p^*(t_k)-p^*(t_{k-1})\|+\|p^*(t_{k-1})-\hat{p}^*\|\leq \frac{2}{3}L,$$
which is a contradiction. \endproof

\subsection{Convergence to the generalized median}
By the discussion in the previous section, we know that the reference points of  all agents that visit an unbounded number of targets converge to a well-defined limit. So do, trivially, the reference points of all agents that visit a bounded set of targets. Hence, we know that the sequence of reference points $p^*_i(t_j)$  converges to a limit $\hat{p}^*_i$, almost surely for all $i \in \{1, \ldots, m\}$. Let use denote by $\mathcal{V}_i(p^*(t_j))$ the Voronoi region associated to the generator $p^*_i(t_j)$, and by  $\mathcal{V}_i(\hat p^*)$ the Voronoi region corresponding to the limit point $\hat{p}^*_i$.
\begin{proposition}
If the limit reference points $\hat p^* = (\hat p_1^*, \ldots, \hat p_2^*)$ are distinct, then the sequence of Voronoi partitions $\{\mathcal{V}(p^*(t_j))\}_{j \in \mathbb{N}}$ converges to the Voronoi partition generated by the limit of reference points, i.e.,
$$\lim_{j \to \infty} \mathcal{V}_i(p^*(t_j)) = \mathcal{V}_i(\hat p^*), \qquad \mathrm{a.s.}$$
\end{proposition}
\proof The boundaries of regions in a Voronoi partition are algebraic curves that depend continuously on the generators, as long as these are distinct. Hence, under this assumption, almost sure convergence of the generators implies the almost sure convergence of the Voronoi regions.
\endproof


As a next step, we wish to understand what is the relation between the asymptotic reference positions and their associated Voronoi regions. More precisely, let $\mathcal{A} \subset \{1, \ldots, m\}$ be the subset of indices of agents that visit an unbounded number of targets; we want to prove that $\hat{p}^*_i$ is indeed the generalized median $\overline{p}_i$ associated to agent $i$, with respect to the limiting set ${\mathcal{V}_j(\hat p^*)}$ and distribution $\varphi(x)$, $\forall i \in \mathcal{A}$. First we need the following technical result.
\begin{lemma}\label{convexity}
Let $\{f_i\}_{i\in \mathbb{N}}:\mathcal{Q}\rightarrow \mathbb{R}$ be a sequence of strictly convex  continuous functions, defined on a common compact subset $\mathcal{Q} \subset \mathbb{R}^n$ . Assume that each $f_i$ has a  unique $x_i:=\arg\min_x f_i$ belonging to the interior of $\Omega$ for any $i$ and that this sequence of function converges uniformly to a continuous strictly convex function $f$ admitting a unique minimum point $\overline{x}$ belonging also to the interior of $\Omega$. Then $\lim_{i\rightarrow \infty}x_i=\overline{x}$.
\end{lemma}
\proof
Since $\{f_i\}_{i\in\mathbb{N}}$ converges uniformly to $f$, then for any $\epsilon>0$, there exists an $I(\epsilon)$ such that for any $i\geq I(\epsilon)$, $\|f_i-f\|\leq \epsilon$ uniformly in $x\in \Omega$. Let $m:=f(\overline{x})$, the minimum value achieved by $f$. Consider the set
$$U_{\epsilon}:=\{x\in \Omega \text{ such that } f(x)\leq m+2\epsilon \}.$$ Since $f$ is strictly convex, for $\epsilon$ sufficiently small, $U_{\epsilon}$ will be a compact subset contained in $\Omega$. Moreover, since  $f$ is strictly convex, we have that $U_{\epsilon}$ is strictly included in $U_{\epsilon '}$, whenever $\epsilon ' > \epsilon$ and both are sufficiently small.
It is also clear that $\lim_{\epsilon\rightarrow \infty} U_{\epsilon}=\overline{x}$ (nested strictly decreasing sequence of compact subsets all containing the point $\overline{x}$). If $\|f_i-f\|\leq \epsilon$, we claim that $x_i\in U_{\epsilon}$; we prove this by contradiction. Since $f_i\geq f-\epsilon$, if $x_i$ does not belong to $U_{\epsilon}$, then $\min(f_i)=f_i(x_i)>m+\epsilon$ (this is just because $U_{\epsilon}$ is simply the set where the function $(f-\epsilon)\leq m+\epsilon$). But since $f_i\leq f+\epsilon$ it turns out that
$f_i(\overline{x})\leq m+\epsilon<\min(f_i)$, which is a contradiction. \endproof

We conclude this section with the following:

\begin{proposition}\label{local}\label{thm:convergence_to_median}
Assume that all agents in $\mathcal{A}$ are able to return infinitely often to their reference point between visiting two targets. Then, the limit reference points of such agents coincide, almost surely, with the generalized medians of their limit Voronoi regions, i.e.,
 $$\hat{p}^*_i= \arg \min_{p \in \Omega} \int_{\mathcal{V}_i(\hat p^*)} \varphi(q) dq, \qquad \mathrm{ a.s. }, \qquad \forall i \in \mathcal{A}.
 $$
\end{proposition}
\proof For any $i\in\mathcal{A}$, define the functions $f_{i,t_j}(p):=\frac{1}{j}\sum_{q\in \mathcal{B}_{i}(t_j)}\|p-q\|$ and $f_i(p):=\int_{\mathcal{V}_i(\hat p^*)}\|p-q\|\varphi(q)\; dq$. These functions are continuous and well defined over $\Omega$. We restrict their domains of definition to the compact set $\mathcal{Q}=\mathrm{Supp}(\varphi)$. These functions are also strictly convex and have unique minima in the interior of $\mathcal{Q}$. Let us notice that, with our previous notation, we have that $p^*_i(t_j)=\arg\min f_{i,t_j}(p)$ and $\overline{p}_i=\arg\min f_i(p)$.
Observe that the functions $f_{i,t_j}(p)$ and $f_i(p)$ can be considered random variables with respect to a probability space whose space of events coincide with all possible realizations of target sequences.
Consider a restriction of these random variable to a new probability space whose space of events coincide with all possible realizations of target sequences, for which the corresponding FT points converge to a limiting point. On this new probability space the random variable $f_i(p)$ becomes a deterministic function which is the expected value of the random variables $f_{i,t_j}(p)$. Since $\mathcal{Q}$ is compact, it is immediate to see that $f_{i,t_j}(p)$ have finite expectation and variance over this reduced probability space, and by the Strong Law of Large Numbers we can conclude that almost surely (over this reduced probability space) $f_{i,t_j}(p)$ converge pointwise to $f_i(p)$. To show that $f_{i,t_j}(p)$ converges pointwise to $f_i(p)$ over the original probability space, it is sufficient to observe the following.
The original probability space is the probability space whose space of events coincide with all possible realizations of target sequences. We already know that almost surely the FT points associated to any possible realization of target sequences will converge. So we can fiber the space of events of the first probability space into spaces of events of reduced probability spaces, except for a set of measure zero.
This is sufficient to prove that $f_{i,t_j}(p)$ converge pointwise to $f_i(p)$ almost surely with respect to all possible realizations of target sequences.

Now that we have proved that almost surely  the sequence $\{ f_{i,t_j}(p)\}_{t_j\in \reals}$ converges pointwise to $f_i(p)$, we prove that it does converge uniformly. To do this, we use a theorem, usually attributed to Dini-Arzela' which state the following: an equicontinuous sequence of  functions converges uniformly to a continuous function on a compact set $\mathcal{Q}$ if and only if it converges poitwise to a continuous function on the same compact set.  Our sequence $f_{i,t_j}(p)$ is equicontinuous if $\forall \epsilon>0$ and $\forall p\in \mathcal{Q}$ there exists a $\delta>0$ such that for all $j\in \{1,\dots,n\}$ and for all $p'\in \mathcal{Q}$ with $\|p'-p\|\leq\delta$, we have $\|f_{i,t_j}(p)-f_{i,t_j}(p')\|\leq\epsilon$; observe that $\delta$ is independent on $j$, while in general it will depend on $\epsilon$ and on $p$.
Now we have $$\|f_{i,t_j}(p)-f_{i,t_j}(p')\|\leq\frac{1}{j}\sum_{q\in \mathcal{B}_{i}(t_j)}\left| \|p-q\|-\|p'-q\|
 \right|.$$
 Using $$\|q-p'\|=\|q-p+p-p'\|\leq \|q-p\|+\|p-p'\| $$
 and $$\|q-p\|=\|q-p'+p'-p\|\leq \|q-p'\|+\|p-p'\|,$$
 it is immediate to see that
 $$\|f_{i,t_j}(p)-f_{i,t_j}(p')\|\leq\frac{1}{j}\sum_{q\in \mathcal{B}_{i}(t_j)}\|p-p'\|\leq \|p-p'\|.$$
 So it is sufficient to take $\delta=\epsilon$ in the previous definition and $\delta$ does not depend on $j$. So the sequence is equicontinuous and the poitwise convergence is upgraded to uniform convergence.


We already know that almost surely the points $p^*_i(t_j)$ do converge to points $\hat{p}^*_i$ (Proposition  (\ref{thm:convergence})); therefore we can claim that $\hat{p}^*_i=\overline{p}_i$, simply applying Lemma (\ref{convexity}), which requires the uniform convergence.
Thus we can claim that the reference position of each agent which services infinitely many targets following our algorithm converges to the generalized median of its Voronoi region, almost surely. \endproof

\subsection{Fairness}
In this section, we prove that, as long as $\varphi$ is strictly positive over a convex set, both policies  introduced in Section \ref{algorithm} are fair.
\begin{proposition}[Fairness]\label{thm:fairness}
If $\mathcal{Q} = \mathrm{Supp}(\varphi)$ is convex, all agents eventually visit an unbounded number of targets, almost surely, i.e.,
$$\lim_{t \to +\infty} \card(\mathcal{B}_i(t)) = +\infty, \qquad \mathrm{a.s.,}\qquad \forall i \in \{1, \ldots, m\}.$$
\end{proposition}
\proof
Under either policy, each agent will pursue and is guaranteed to service
the nearest target in its own Voronoi region. Hence, in order to show that an agent services an unbounded number of targets, it is sufficient to show that the probability that the next target be generated within its Voronoi region remains
strictly positive, i.e., that
$$\int_{\mathcal{V}_i(p(t))} \varphi(q) dq > 0.$$
Since $\mathcal{Q}$ is convex, and all agents move towards the nearest target, at least initially,  all agents will eventually enter $\mathcal{Q}$, and remain within it. Let us denote with $P_{ij}$ the probability that agent $i$ visits the  $j$-th target. For any $i$ this probability is always strictly positive. Indeed, even if it happens that some agents are servicing simultaneously the same target (simply because the service request appears at the boundary of two different Voronoi regions), this does not mean that their reference points have to coincide. For the reference points of two agents eventually to coincide, it must happen that they are servicing infinitely many often  and simultaneously the same target. Since the boundaries of Voronoi regions have measure zero  and $\varphi$ is a continuous distribution without singular components, we can claim that  $\lim_{j\rightarrow \infty}P_{ij}>0$ almost surely.
Now call $P=\min_{i=1,\dots,m}\lim_{j\rightarrow\infty}P_{ij}$. Then $P$ is strictly positive almost surely.
Therefore, the probability that the $i$-th agent {\em does not} visit an unbounded number of targets is bounded from above by $\lim_{j\rightarrow\infty}\prod_{k=1}^j (1-P)^k=0$, almost surely. \endproof

\subsection{Efficiency}
In this section, we will prove that the system time provided by either one of the algorithms in Section \ref{algorithm} converges to a critical point (either a saddle point or  a local minimum) with high probability as $\lambda \to 0$.


In the preceding sections, we have proved that---as long as each agent is able to return to its reference point between servicing two targets, infinitely often---the reference points $p^*_i(t_j)$ converge to points $\hat{p}^*_i$, which generate a MVT. In such case, we know that the average time of service will converge to
\begin{equation}\label{eq:atos}
\overline T_\pi =  \int_\Omega \min_{i=1, \ldots, m} \left\|  \hat{p}^*_i - q\right\| \; \varphi(q) dq  = \sum_{i=1}^m \int_{\mathcal{V}_i(\hat{p}^*)} \|\hat{p}^*_i - q\| \; \varphi(q) dq.
\end{equation}
Consider now  functions $\mathcal{H}_m$ of the form:
\begin{equation}\label{eq:Hm}
\mathcal{H}_m(p_1,\dots ,p_m)= \int_\Omega \min_{i=1, \ldots, m} \left\|  p_i - q\right\| \; \varphi(q) dq  = \sum_{i=1}^m \int_{\mathcal{V}_i(p)} \|p_i - q\| \; \varphi(q) dq.
\end{equation}
Observe that $\overline{T}_\pi$ belongs to the class of functions of the form $\mathcal{H}_m$ where each point $p_i$ is constrained to be the generalized median of the corresponding Voronoi region (i.e., $\mathcal{V}(p)$ is a MVT).

We want to prove that $\overline{T}_\pi$ is a critical point of ${\cal H}_m$. To do so, we consider an extension of ${\cal H}_m$, i.e. a functional ${\cal K}_m$ defined as follows:
$${\cal K}_m(p_1,\dots,p_m,\mathcal{V}_1,\dots,\mathcal{V}_m):=\sum_{i=1}^m\int_{y\in \mathcal{V}_i}\left\|y-p_i\right\| \; \varphi(y)dy.$$
Observe that in this case the regions $\{\mathcal{V}_i\}_{i=1,\dots,m}$ are not restricted to form a MVT with respect to the generators $\{x_i\}_{i=1,\dots,m}$. Thus we can view the functional ${\cal H}_m$ we are interested in as a constrained form of the unconstrained functional ${\cal K}_m$. It turns out therefore that critical points of ${\cal K}_m$ are also critical points of ${\cal H}_m$. With respect to critical points of ${\cal K}_m$ we have the following result:

\begin{proposition}\label{efficiency}\label{thm:critical}
Let $\{p_i\}_{i=1,\dots,m}$ denote {\em any set} of $m$ points belonging to $\mathrm{Supp}(\varphi)$ and let $\{\mathcal{V}_i\}_{i=1,\dots,m}$ denote {\em any tessellation} of $\mathrm{Supp}(\varphi)$ into $m$ regions. Moreover, let us define ${\cal K}_m$ as above.

Then a sufficient condition for $\{p_1,\dots,p_m,\mathcal{V}_1,\dots,\mathcal{V}_m\}$ to be a critical point (either a saddle point or a local minimum), is that the $\mathcal{V}_i$'s are the Voronoi regions corresponding to the $p_i$'s, and, {\em simultaneously}, the $p_i$'s are the generalized median of the corresponding $\mathcal{V}_i$'s.
\end{proposition}
\proof Consider first the variation of ${\cal K}_m$ with respect to a single point, say $p_i$. Now let $v$ be a vector in $\mathbb{R}^2$, such that $p_i+\epsilon v \in \Omega$. Then we have
$$ {\cal K}_m(p_i+\epsilon v)-{\cal K}_m(p_i)=\int_{y\in \mathcal{V}_i}\left\{\|y-p_i-\epsilon v\|-\|y-p_i\|\right\}\;\varphi(y) dy,$$
where we have not listed the other variables on which ${\cal K}_m$ depends since they remain constant in this variation. By the very form of this variation,  it is clear that if the point $p_i$ is the generalized median for the {\em fixed} region $\mathcal{V}_i$, we will have that  ${\cal K}_m(p_i+\epsilon v)-{\cal K}_m(p_i)>0$, for any $v$.
Now consider the points $\{p_i\}_{i=1,\dots,m}$ fixed and consider a tessellation $\{\mathcal{U}_i\}_{i=1,\dots,m}$ different from the Voronoi regions $\{\mathcal{V}_i\}_{i=1,\dots,m}$ generated by the points $p_i$'s. We compare the value of ${\cal K}_m(p_1,\dots,p_m,\mathcal{V}_1,\dots,\mathcal{V}_m)$, with the value of ${\cal K}_m(p_1,\dots,p_m,\mathcal{U}_1,\dots,\mathcal{U}_m).$ Consider those $y$ which belong to the Voronoi region $\mathcal{V}_j$ generated by $p_j$, and possibly not to the Voronoi region of another $p_i$. Anyway, since $\mathcal{U}_i$ is not a Voronoi tessellation, it can happen that in any case these $y$ belong to $\mathcal{U}_i$. Thus for these particular $y$'s we have $\varphi(y)\|y-p_j\|\leq\varphi(y)\|y-p_i\|$. Moreover, since $\{\mathcal{U}_i\}_{i=1,\dots,m}$ are not the Voronoi tessellation associated to the $p_i$'s, the last inequality must be strict over some set of positive measure. Thus we have that ${\cal K}_m(p_1,\dots,p_m,\mathcal{V}_1,\dots,\mathcal{V}_m)<
{\cal K}_m(p_1,\dots,p_m,\mathcal{U}_1,\dots,\mathcal{U}_m)$, and therefore ${\cal K}_m$ is minimized, keeping fixed the $p_i$'s exactly when the subset $\mathcal{V}_i$'s are chosen to be the Voronoi regions associated with the point $p_i$'s. \endproof

By the previous proposition and by the fact that critical points of the unconstrained functional ${\cal K}_m$ are also critical points of the constrained functional ${\cal H}_m$, we have that the MVT are always critical points for the functional ${\cal H}_m$, and in particular $\overline{T}$ is either a saddle point or a local minimum for the functional ${\cal H}_m$.

Before we conclude, we need one last intermediate result.
\begin{proposition}
Each agent will be able to return to its reference point before the generation of a new service request infinitely often with high probability as $\lambda \to 0$.
\end{proposition}
\proof
Let $t_1$ be such that $\mathcal{B}_i(t_1) \neq \emptyset$, for all $i \in \{1, \ldots, m\}$. Such time exists, almost surely, because of the fairness of the proposed policies. At time $t_1$ all agents will be within $\mathcal{Q}$. Let $n_1 = \mathrm{card}(\mathcal{D}(t_1))$ be the total number of outstanding targets at time $t_1$. An upper bound on the time needed
to visit all targets in $\mathcal{D}(t_1)$ is  $n_1 (\diam(\mathcal{Q})$.
When there are no outstanding targets, agents move to their reference points, reaching them in at most $\diam(\mathcal{Q})$ units of time.

The time needed to service the initial targets and go to the reference configuration is hence bounded by $t_\mathrm{ini}\le t_1 + (n_1+1) \diam(\mathcal{Q})$.  The probability that at the end of this initial
phase the number of targets is reduced to zero is
$$  P\left[n(t_\mathrm{ini})=0\right] = \exp(-\lambda (t_\mathrm{ini}-t_1))
  \geq \exp(-\lambda (n_0+1) \diam(\mathcal{Q})),$$
that is, $P\left[n(t_\mathrm{ini})=0\right] \rightarrow 1^-$ as $\lambda
\rightarrow 0^+$. As a consequence, after an initial transient, all targets
will be generated with all agents waiting at their reference points, and an empty demand queue, with high probability as $\lambda \to 0^+$.
\endproof

We can now conclude with following:

\begin{theorem}[Efficiency]
The system time provided by the no-communication policy $\pi_\mathrm{nc}$ and by the sensor-based policy $\pi_\mathrm{sb}$ converges to a critical point (either a saddle point or a local minimum) with high probability as $\lambda \to 0$.
\end{theorem}
\proof Combining results in Propositions \ref{thm:convergence} and \ref{thm:convergence_to_median} we conclude that the reference points of all agents that visit an unbounded number of targets converge to a MVT, almost surely---provided agents can return to the reference point between visiting targets. Moreover, the fairness result in Proposition \ref{thm:fairness} shows that in fact all agents do visit an unbounded number of targets almost surely; as a consequence, Proposition \ref{thm:critical} the limit configuration is indeed a critical point for the system time.  Since agents return infinitely often to their reference positions with high probability as $\lambda \to 0$, the claim is proven.
\endproof

Thus we have proved that the suggested algorithm enables the agents to realize a coordinated task, such that ``minimizing" the cost function without explicit communication, or with mutual position knowledge only. Let us underline that, in general, the achieved critical point strictly depends on the initial positions of the agents inside the environment $\Omega$. It is known that the function $\mathcal{H}_m$ admits (not unique, in general) {\em global minima}, but the problem to find them is {\em NP-hard}.

\begin{remark}We can not exclude that the algorithm so designed will converge indeed to a saddle point instead of a local minimum. This is due to the fact that the algorithm provides a sort of implementation of the {\em steepest descent} method, where, unfortunately we are not following the steepest direction of the gradient of the function $\mathcal{H}_m$, but just the gradient with respect to one of the variables.
For a broader point of view of steepest descent in this framework see for instance  \cite{Okabe.Boots.ea:00}.
\end{remark}

On the other hand, since the algorithm is based on a sequence of targets and at each phase we are trying to minimize a different cost function, it can be proved that the critical points reached by this algorithm are {\em no worse} than the critical points reached knowing a priori the distribution $\varphi$. This is a remarkable result proved in a different context in \cite{Sabin.Gray:86}, where it is also presented an example in which the use of a sample sequence provides a better result (with probability one) than the a priori knowledge of $\varphi$. In that specific example the algorithm with the sample sequence does converge to a global minimum, while the algorithm based on the a priori knowledge of the distribution $\varphi$ gets stuck in a saddle point.

\subsection{A comparison with algorithms for vector quantization and centroidal Voronoi tessellations}

The use of Voronoi tessellations  is ubiquitous in many fields of science, ranging from operative research, animal ethology (territorial behaviour of animals),  computer science (design of algorithms), to numerical analysis (construction of adaptive grids for PDEs and general quadrature rules), and algebraic geometry (moduli spaces of abelian varieties). For a detailed account of possible applications see for instance the book \cite{Okabe.Boots.ea:00}.
In the available literature, most of the analysis is devoted to applications of centroidal Voronoi tessellations, i.e., Voronoi tessellation such that
$$p_i = \arg \min_{s \in \mathbb{R}^n} \int_{\mathcal{V}_i(p_1, p_2, \ldots, p_m)} \|s-q\|^2 \varphi(q) \; dq, \qquad \forall i \in \{1, \ldots, m\}.$$

A popular algorithm due to Lloyd~\cite{Lloyd:82} is based on the iterative computation of centroidal Voronoi tessellations. The algorithm can be summarized as follows. Pick $m$ generator points, and consider a large number $n$ of samples from a certain distribution. At each step of the algorithm generators are moved towards the centroid of the samples inside their respective Voronoi region. The algorithm terminates when each generator is within a given tolerance from the centroid of samples in its region, thus obtaining a centroidal Voronoi tessellation weighted by the sample distribution.  There is also a continuous version of the algorithm, which requires the a priori knowledge of a spatial density function, and computation of the gradient of the polar moments of the Voronoi regions with respect to the positions of the generators. An application to coverage problems in robotics and sensor networks of Lloyd's algorithm is available in ~\cite{Cortes.Martinez.ea:04}.

The algorithms we introduced in this paper are more closely related to an algorithm due to MacQueen~\cite{MacQueen:67}, originally designed as a simple online adaptive algorithm to solve clustering problems, and later used as the method of choice in several vector quantization problems where little information about the underlying geometric structure is available.
MacQueen's algorithm can be summarized as follows. Pick $m$ generator points. Then iteratively sample points according to the probability density function $\varphi$. Assign the sampled point to the nearest generator, and update the latter by moving it in the direction of the sample. In other words, let $q_j$ be te $j$-th sample, let ${i^*(j)}$ be the index of the nearest generator, and let $c = (c_1, \ldots, c_m)$ be a counter vector, initialized to a vector of ones. The update rule takes the form
$$p_{i^*(j)} \leftarrow \frac{c_{i^*(j)} q_j + p_{i^*(j)}}{c_{i^*(j)}+1},$$
$$c_{i^*(j)} \leftarrow c_{i^*(j)}+1.$$
The process is iterated until some termination criterion is reached.  Compared to Lloyd's algorithm,  MacQueen's algorithm has the advantage to be a learning adaptive algorithm, not requiring the a priori knowledge of the distribution of the objects, but rather allowing the online generation of samples. It can be recognized that the update rule in MacQueen's algorithm  corresponds to moving the generator points to the centroids of the samples assigned to them.

The algorithm we propose is very similar in spirit to MacQueen's algorithm, however, there is a remarkable difference. MacQueen's algorithm deals with centroidal Voronoi tessellations, thus with the computation of $k$-means. Our algorithm instead is based on MVT, and on the computation of $k$-medians. In general, very little is known for Voronoi diagrams generated using simply the Euclidean distance instead of the square of the Euclidean distance.
For instance, the medians of a sequence of points can exhibit a quite peculiar behavior if compared to the one of the means. Consider the following example.
Given a sequence of points $\{q_i\}_{i\in \mathbb{N}}$ in a compact set $K\subset \mathbb{R}^2$, we can construct the induced sequence of means:
$$ m_N:=\frac{1}{N}\sum_{i=1}^N q_j$$
and analogously the induce sequence of FT points we considered in the previous sections. Call $FT_N$ the FT point corresponding to the first $N$ points of the sequence $\{q_i\}_{i\in\mathbb{N}}$.
We want to point out that induced sequence $\{m_j\}$ and $\{FT_j\}$ have a very different behaviour. Indeed, the induced sequence of means will always converge as long as the points $q_j$s belong to a compact set. To see this, just observe that if $\mathrm{diam}(K)\leq L$, then $\|m_j\|\leq L$. Moreover,
it is immediate to see that $\|m_{N+1}-m_{N}\|\leq \frac{2L}{N}$. Then one can conclude using the same argument of Theorem (\ref{thm:convergence}).
On the other hand, one can construct a special sequence of points $q_j$s in a compact set $K$ for which the induced sequence of FT points does not converge. This is essentially due to the fact that while the contribution of each single point $q_j$ in moving the position of the mean decreases as $j$ increases, this could not happen in the case of the median. To give a simple example, start with the following configuration of points in $\mathbb{R}^2$: $q_1=(1,0), q_2=(-1,0), q_3=(0,1)$ and $q_4=(0,-1)$. Then the sequence of points $q_j$s continue in the following way: $q_k=(0,1)$ if $k>4$, and $k$ odd, $q_k=(0,-1)$ if $k>4$ and $k$ is even. Using the characterization of FT points, it is immediate to see that $FT_k=(0,0)$ for $k>4$ and $k$ even, while $FT_k=(0,\tan(\pi/6))$ for $k>4$ and $k$ odd, so the induced sequence can not converge. This phenomenon can not happen to the sequence of means, which is instead always convergent. Therefore, it should be clear that the use of MVT instead of centroidal Voronoi tessellations makes much more difficult to deal with the technical aspects of the algorithm such as its convergence.

\section{A game-theoretic point of view}\label{games}
In this section we provide an analysis of the proposed algorithm from the point of view of game theory. In particular, we frame our presentation on the works \cite{Arslan.Shamma:06}, \cite{Shamma.Arslan:05}, in which the point of view of game theory has been introduced in the study of cooperative control and strategic coordination of decentralized networks of multi-agents systems.
In this section we prove that our algorithm provides a pure Nash equilibrium in a multi-player game where each agent is interested in maximizing its own utility. On the other hand, our multi-player game formulation is much simpler than the usual framework considered in the literature, since there will be no negotiation mechanism among the agents. Despite this fact it turns out that in this example, just trying to maximize their own utility function the agents will indeed maximize a {\em different global utility function}.

In this section we view the agents as {\em rational autonomous decision makers} trying to maximize their own utility function. The utility function of agent $i$, denoted by $U_i$ is simply the {\em expected} number of service requests handled by agent $i$ within a certain time horizon, bounded or unbounded, where $i=1, \dots, m$. We assume that the stochastic process for generating targets is the one already described in the previous sections.
It is obvious that any of the utility function $U_i$ is a function of the policy vector $\pi:=\{\pi_1,\dots,\pi_m\}$ consisting of all the policies followed by each agent. In general the space of all policies $\Pi$ is just an uncountable set containing all conceivable policies chosen by an agent and it does not have any other structure. Thus we assume that the policy space $\Pi_i$ of agent $i$ is equal to a fixed policy space $\Pi$ which is independent on $i$, this for any $i=1,\dots, m$. Therefore the policy vector $\pi
\in\Pi^{m}$. Denoting with $\pi_{-j}:$ the policy specification of all the agents, except agent $j$, i.e.
$\pi_{-j}:=(\pi_1,\dots,\pi_{j-1},\pi_{j+1},\dots,\pi_m)$ we may write policy vector $\pi$ as $(\pi_j,\pi_{-j})$.
Using this notation we can formulate the following definition adapted by \cite{Arslan.Shamma:06}:
\begin{definition}\label{pureNash}
A policy vector $\pi^*$ is called a {\em pure Nash equilibrium}  if for all $j=1,\dots,m$:
\begin{equation}\label{pureeq}
U_j(\pi_j^*,\pi_{-j}^*)=\max_{\pi_j \in \Pi}U_j(\pi_j,\pi_{-j}^*).
\end{equation}
Moreover, a policy vector $\pi$ is called {\em efficient} if there is no other policy vector that yields higher utilities to all agents.
\end{definition}

Under the target generation assumptions followed so far, we have the following:
\begin{proposition}
Let us call $\tilde{\pi_j}$ the policy assignment for agent $j$ corresponding to our algorithm. Then the policy vector $\tilde{\pi}:=\{\tilde{\pi_j}\}_{j=1,\dots,m}$ is an efficient pure Nash equilibrium for the given agents utilities.
\end{proposition}
\proof It is immediate to see that policy vector $\tilde{\pi}$ satisfies equation (\ref{pureeq}), and thus is a pure Nash equilibrium. It is also clear that cannot be any other policy assignment which yields {\em strictly higher} utilities to all agents, simply because when there is a new outstanding service request all agents move directly toward that location, trying to satisfy it as if there were no other agent in the environment. \endproof

 On the other hand, observe that we do {\em not} claim in the previous proposition that our algorithm provides  the unique efficient pure Nash equilibrium for the given set of utilities function. For instance to find other efficient pure Nash equilibria it is sufficient to modify  the algorithm during the initial phases, when the FT points are not uniquely determined. These modifications produce different policy vectors which are anyway all efficient pure Nash equilibria, as it is immediate to see.

 Our game-theoretic formulation of the given algorithm belongs to a class of multi-player games called potential games. In a potential game the difference in the utility reached by any of the agents for two different policy choices, when the policies of the other agents are kept fixed, can be measured by a potential function that depends on the policy vector and not on the label of any agent. The fact that our game-theoretic formulation belongs to this class is obvious since the functional form of the utility function is the same for each agent, that is $U_j=U$ for any agent $j$. So in this case, as a potential function we take $U$. The formal definition is as follows:

 \begin{definition}\label{PG}
 A {\em potential game} is the set $\{U_1,\dots,U_m,\psi\}$ consisting of agent utilities $U_1(\pi),\dots,U_n(\pi)$ and a potential function $\psi:\Pi^m\rightarrow \mathbb{R}$, such that for every agent $a_j$ and for any policy assignments $\pi_j', \pi_j"\in \Pi_j$ and $\pi_{-j}\in \prod_{k\neq j}\Pi_k$:
 $$U_j(\pi_j^{'},\pi_{-j})-U_j(\pi_j^{"},\pi_{-j})=\psi(\pi_j^{'},\pi_{-j})-\psi(\pi_j^{"},\pi_{-j}).$$
 \end{definition}

 An extension of this concept is provided by:
 \begin{definition}\label{OPG}
 An {\em ordinal potential game} is the set $\{U_1,\dots,U_m,\psi\}$ consisting of agent utilities $U_1(\pi),\dots,U_n(\pi)$ and a potential function $\psi:\Pi^m\rightarrow \mathbb{R}$, such that for every agent $a_j$ and for any policy assignments $\pi_j', \pi_j"\in \Pi_j$ and $\pi_{-j}\in \prod_{k\neq j}\Pi_k$:
 $$U_j(\pi_j^{'},\pi_{-j})-U_j(\pi_j^{"},\pi_{-j})>0 \text{ if and only if } \psi(\pi_j^{'},\pi_{-j})-\psi(\pi_j^{"},\pi_{-j})>0.$$
 \end{definition}

 At this point, we are ready to introduce the global utility function for this multi-player game formulation.
 The global utility function for this game is given by $U_g(\pi)=-\overline{T}_\pi$, where $\overline{T}_\pi$ is the {\em system time} under policy vector $\pi$
 An important aspect in the game-theoretic approach is to understand to which extent the utility functions of the individual players are compatible with the global utility function. To this aim the following definition has been introduced in \cite{Arslan.Shamma:06}:

 \begin{definition}\label{aligned}
 The set of agents utilities $\{U_j(\pi)\}_{j=1,\dots,m}$ is {\em aligned} with the global utility $U_g(\pi)$ iff
 the set $\{U_1,\dots,U_m,U_g\}$ forms an ordinal potential game, with potential function given by $U_g$.
 \end{definition}

 It is immediate to prove the following:

 \begin{proposition}
The given agents utilities are aligned with the global utility function, in this game-theoretic formulation of the algorithm.
 \end{proposition}
 \proof Let us focus on one agent, say agent $j$. If its utility function increases, it is able to service a bigger number of targets in the given time horizon. It can happen that some other agent will have a corresponding decrease in their utility function, and this happens exactly when the agent $j$ due to a policy change will service some targets more rapidly. This in turn, will increase anyway $U_g$. \endproof

In general, it is true that alignment does not prevent pure Nash equilibria from being suboptimal from the point of view of the global utility. Moreover, even {\em efficient} pure Nash equilibria (i.e. pure Nash equilibria which yield the highest utility to all agents) can be suboptimal from the perspective of the global utility function. Such a phenomenon is indeed what happens in our construction.

\begin{proposition}
The policy vector $\tilde{\pi}$ corresponding to the policies realized by our algorithm is an efficient pure Nash equilibrium which is possibly suboptimal from the point of view of global utility.
\end{proposition}
\proof We already know that $\tilde{\pi}$ is an efficient pure Nash equilibrium and from the analysis developed in Section (\ref{prova}) we know that it yields a {\em critical point} for the system time $\overline{T}_\pi$. Thus it corresponds to a {\em critical point} for $U_g(\pi)$, either a local maximum or a saddle point. On the other hand, as noted for  $\overline{T}_\pi$, $U_g$ may have in general several local maxima and our algorithm is not guaranteed to converge to a global maximum. \endproof

\section{Numerical results}
\label{simulations}
In this section, we present simulation results showing the performance of the proposed policies for various scenarios.

\subsection{Uniform distribution, light load}
In the numerical experiments, we first consider $m=9$, choose $\mathcal{Q}$ as a unit square, and set $\varphi = 1$ (i.e., we consider a spatially uniform target-generation process). This choice allows us to determine easily the optimal placement of reference points, at the centers of a tesselation of $\mathcal{Q}$ into nine equal squares, and compute analytically the optimal system time.
In fact, it is known that the expected distance of a point $q$ randomly sampled from a uniform distribution within a square of side $L$ from the center of the square $c$ is
$$\mathrm{E}[\|q- c\|] = \frac{\sqrt{2} + \log(1 + \sqrt{2})}{6} L \approx 0.3826 L.$$

\begin{figure}[tbh]
\null \hfill
\includegraphics[width=0.4\textwidth]{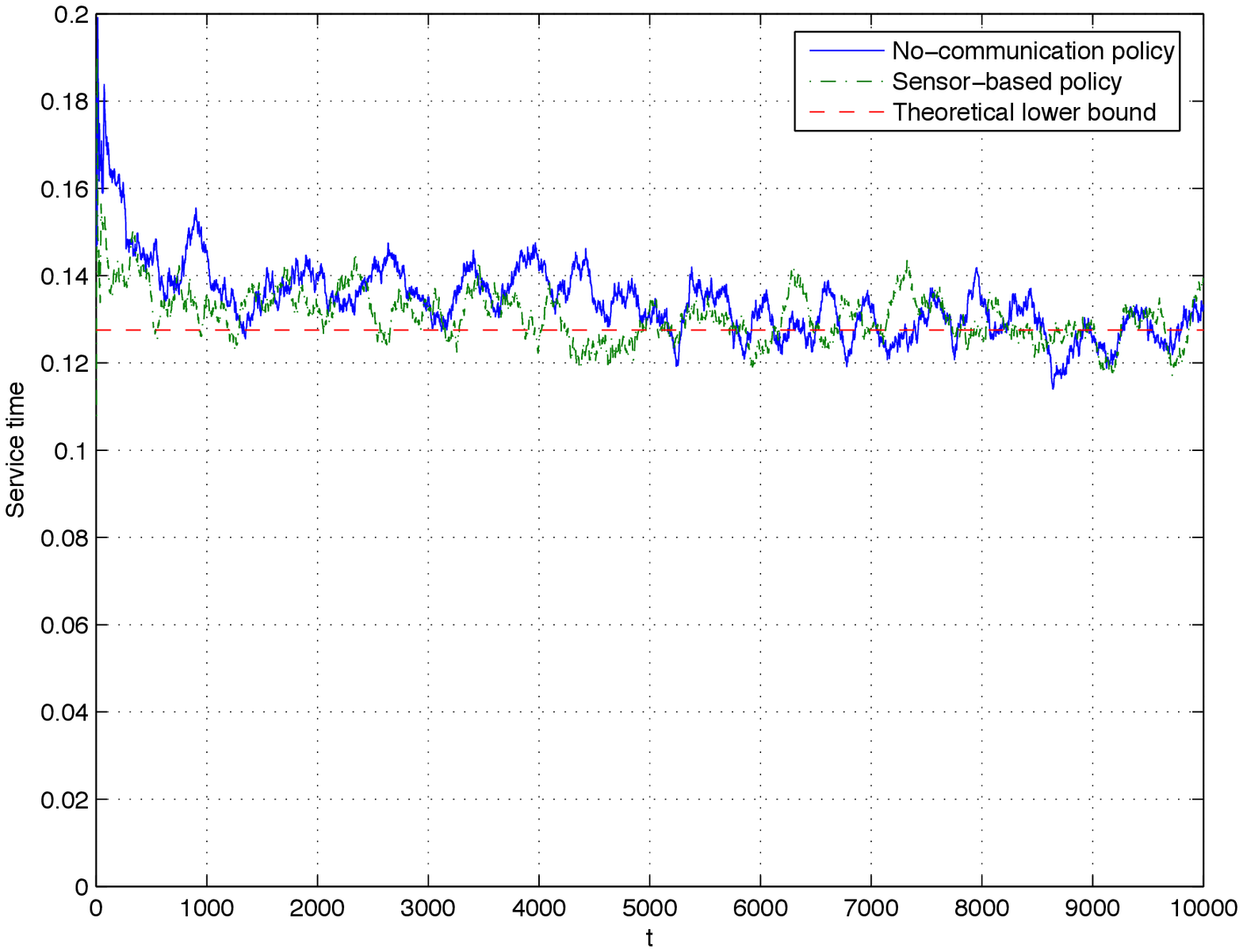}
\hfill \includegraphics[width=0.4\textwidth]{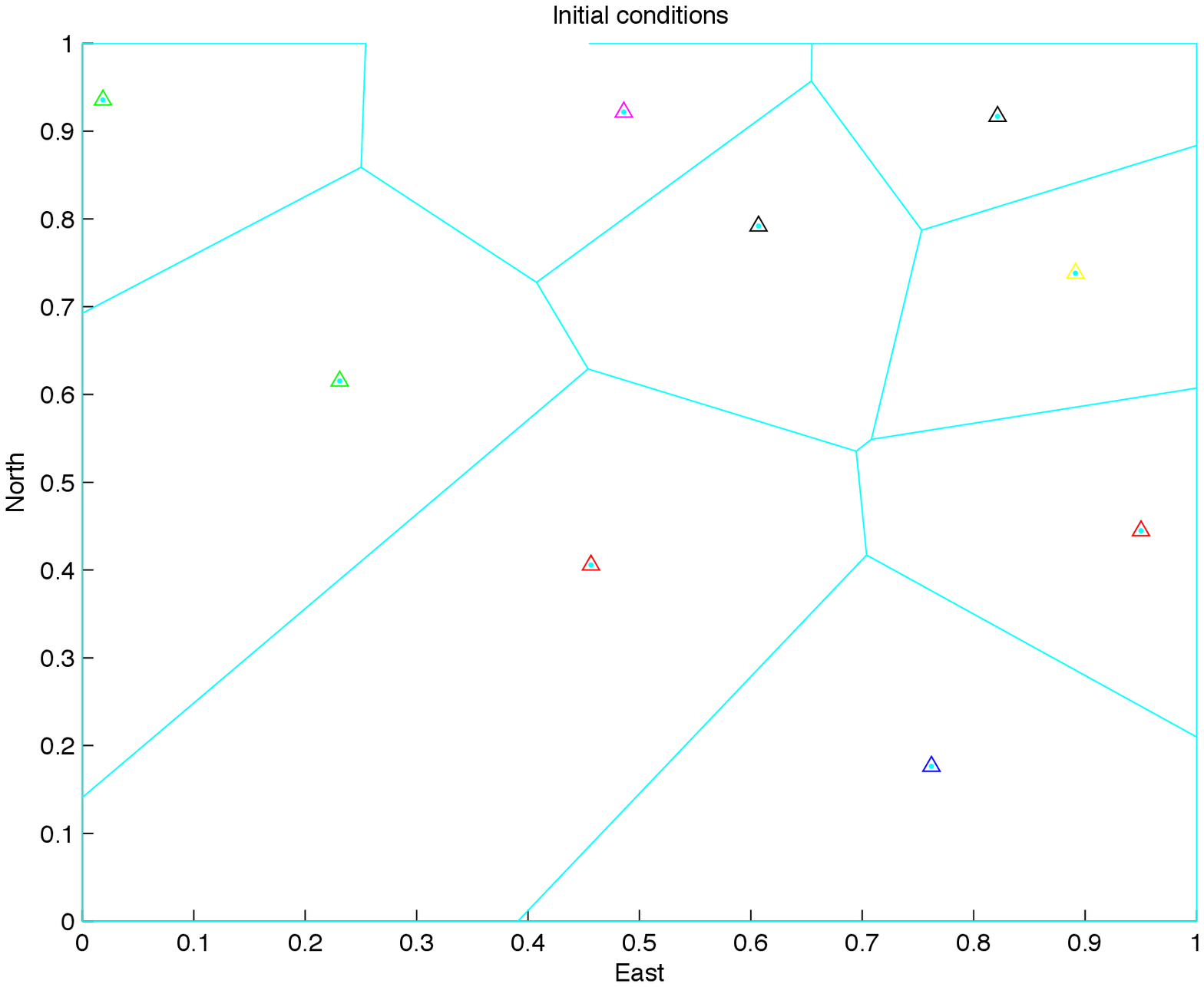} \hfill \null\\
\null\hfill \includegraphics[width=0.4\textwidth]{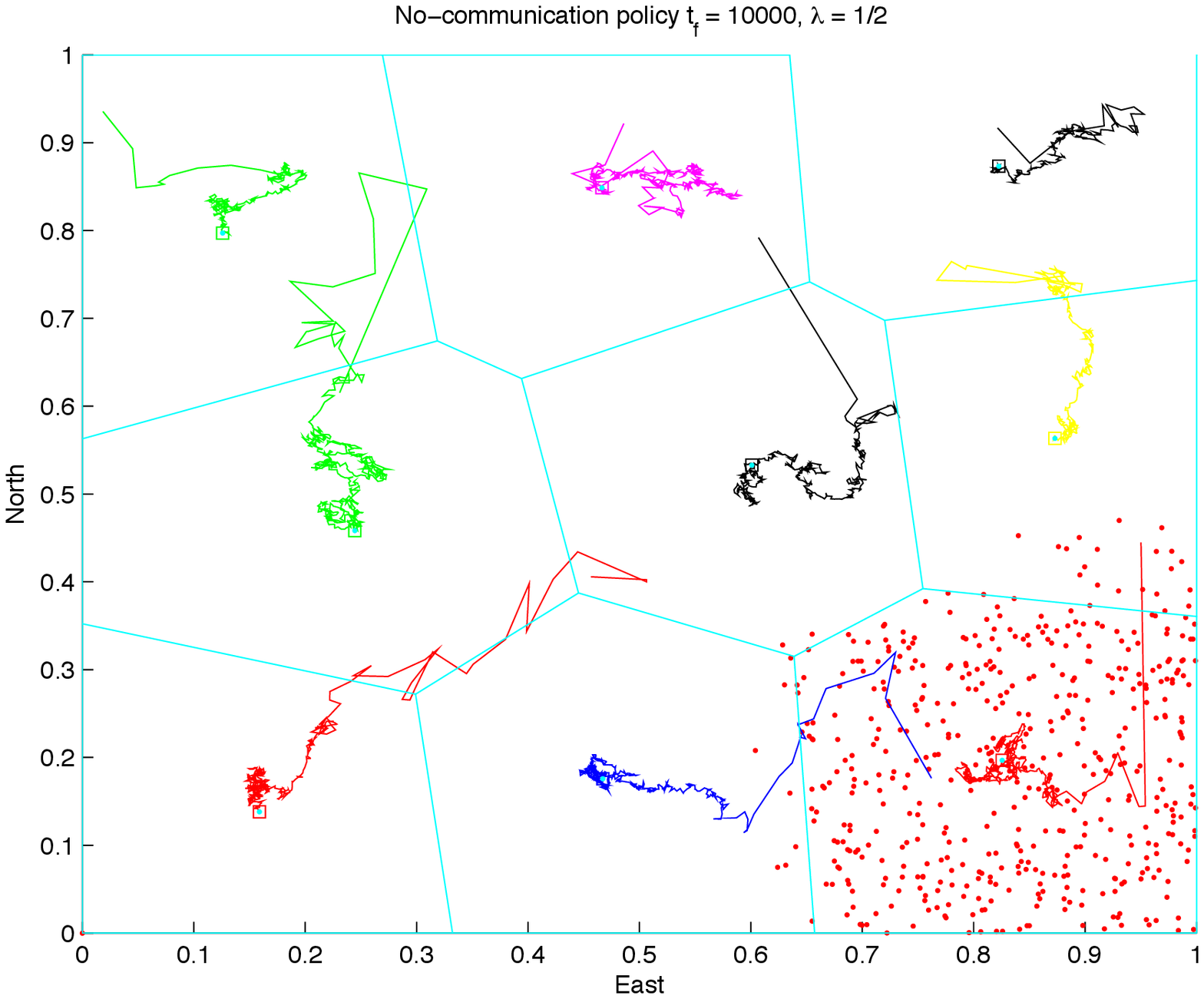}
\hfill \includegraphics[width=0.4\textwidth]{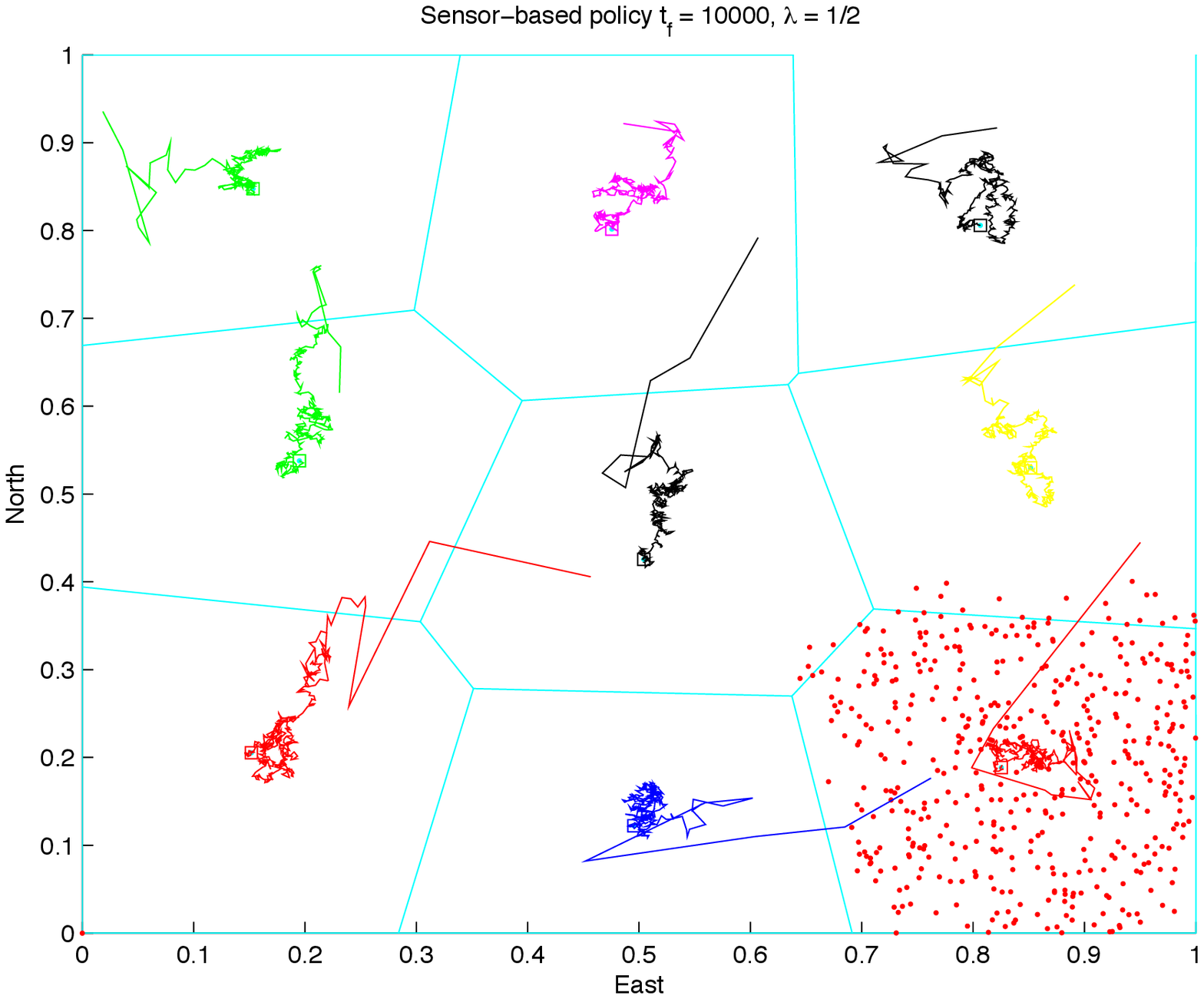} \hfill\null
\caption{Numerical simulation in the light-load case, for a uniform spatial distribution. Top left: the actual service times as a function of time, for the two policies, compared with the optimal system time. Top right: the initial configuration of the nine agents. Bottom left and right: paths followed by the reference points up to $t=10^4$ (corresponding to approximately 5,000 targets), using the two policies. The locations of all targets visited by one of the agents are also shown.}
\label{fig:light_load}
\end{figure}

The results for a small value of $\lambda$, i.e., $\lambda = 0.5$, are presented in Figure \ref{fig:light_load}.
The average service time converges to a value that is very close to the theoretical limit computed above, taking $L  = 1/\sqrt{m} = 1/3$. In both cases, the reference points converge---albeit very slowly---to the generators of a MVT, while the average system time quickly approaches the optimal value.

\subsection{Non-uniform distribution, light load}
\begin{figure}[tbh]
\null \hfill
\includegraphics[width=0.4\textwidth]{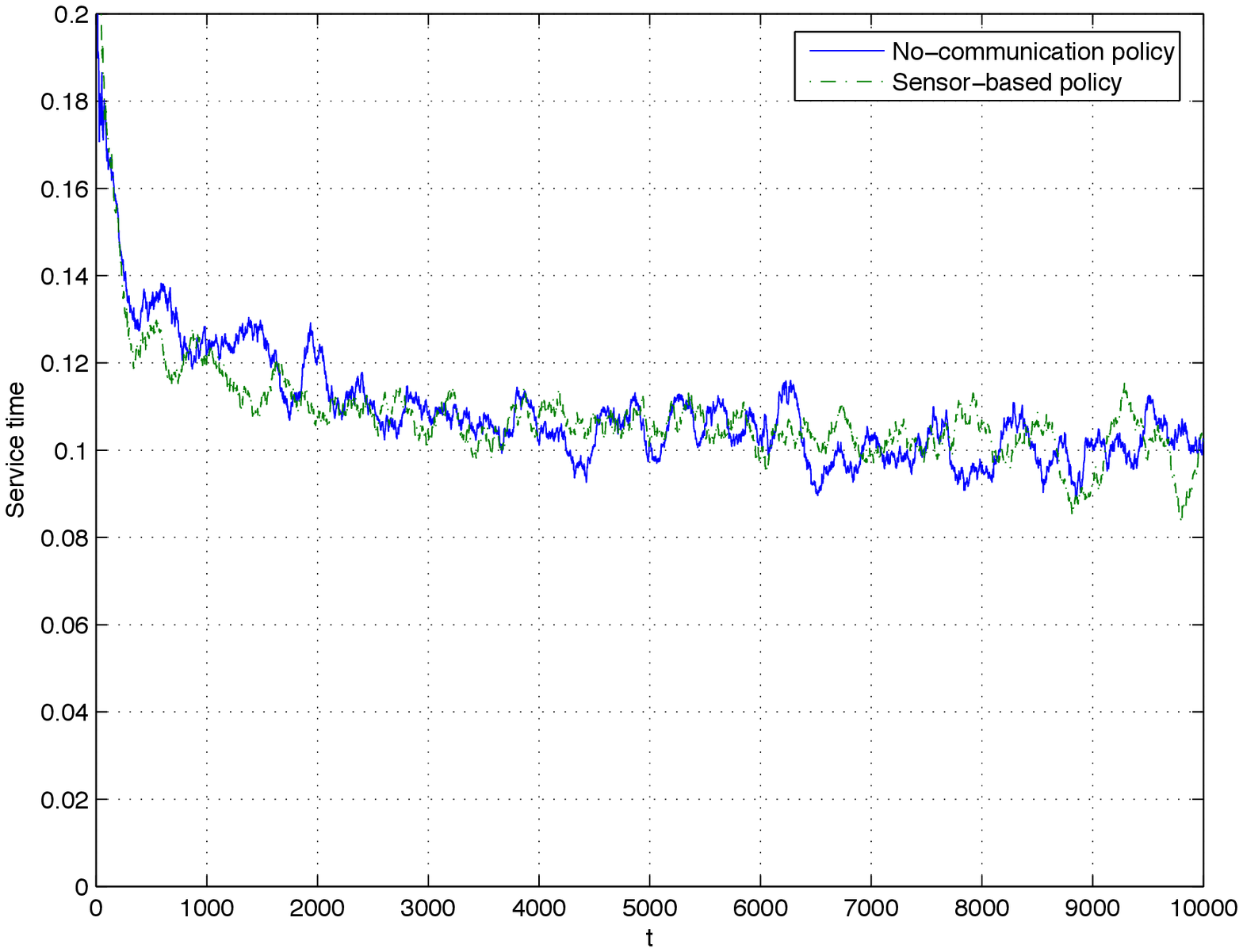}
\hfill \includegraphics[width=0.4\textwidth]{refpaths-ini} \hfill \null\\
\null\hfill \includegraphics[width=0.4\textwidth]{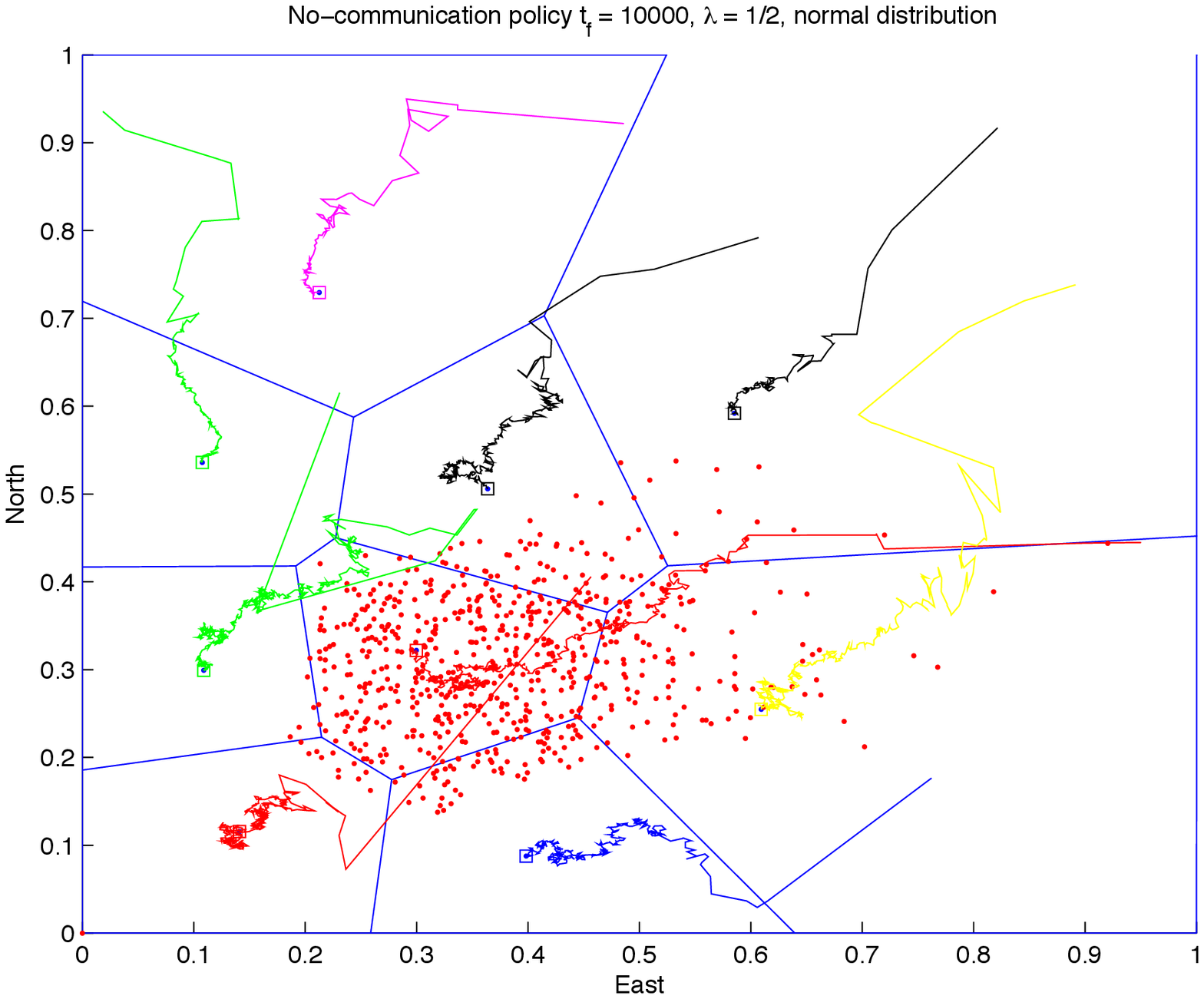}
\hfill \includegraphics[width=0.4\textwidth]{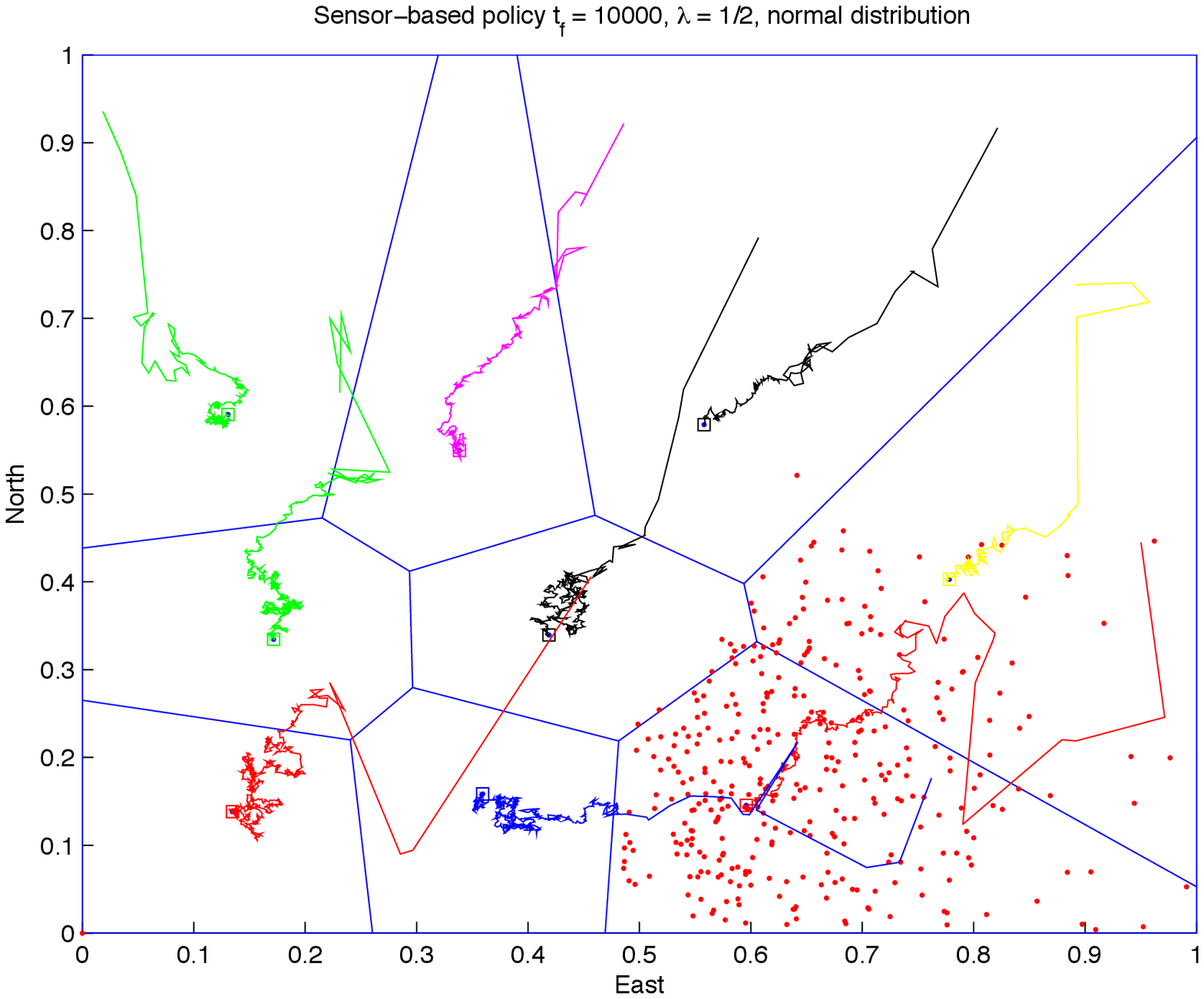} \hfill\null
\caption{Numerical simulation in the light-load case, for a normal spatial distribution. Top left: the actual service times as a function of time, for the two policies. Top right: the initial configuration of the nine agents. Bottom left and right: paths followed by the reference points up to $t=10^4$ (corresponding to approximately 5,000 targets), using the two policies. The locations of all targets visited by one of the agents are also shown.}
\label{fig:light_load_normal}
\end{figure}

We also present in Figure \ref{fig:light_load_normal} results of similar numerical experiments with a non-uniform distribution, namely an isotropic normal distribution centered at $(0.25,0.25)$, with standard deviation equal to $0.25$.

\subsection{Uniform distribution, dependency on the target generation rate}
\begin{figure}[htb]
\centerline{\includegraphics[width=0.8\textwidth]{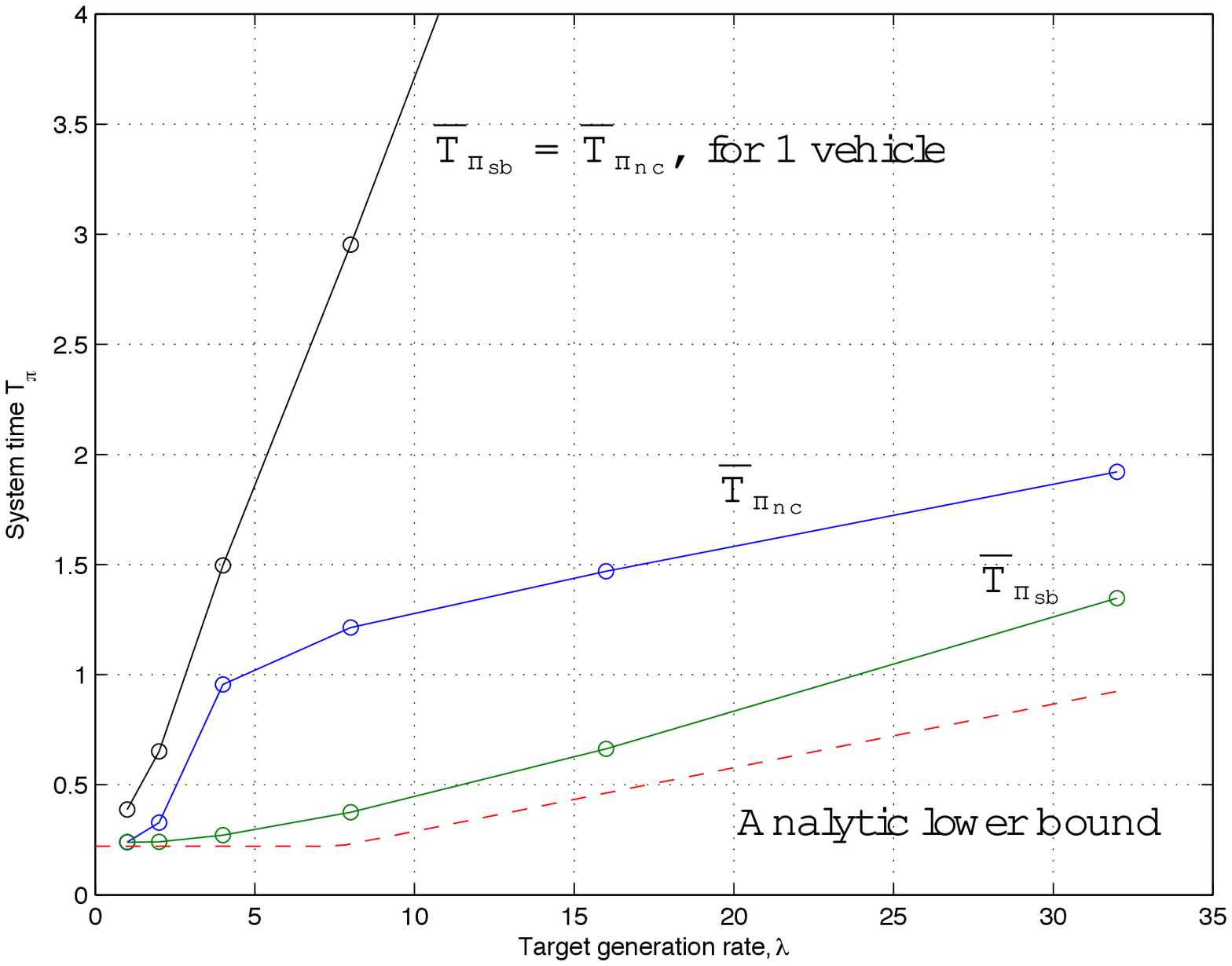}}
\caption{System time provided by the policies proposed in this
paper, as a function of the target generation rate $\lambda$. The
system is composed of three vehicles, and the target points are
generated uniformly in the unit square.} \label{fig:TvsLambda}
\end{figure}

An interesting set of numerical experiments evaluates the performance of the proposed policies over a large range of values of the target generation rate $\lambda$. In Section \ref{prova}, we proved the convergence of the system's behavior to an efficient steady state, with high probability as $\lambda \to 0$, as confirmed by the simulations discussed above. For large values of $\lambda$ however,  the assumption that vehicles are able to return to their reference point breaks down, and the convergence result is no longer valid. In figure \ref{fig:TvsLambda} we report results from numerical experiments on scenarios involving $m=3$ agents, and values of $\lambda$ ranging from $1/2$ to $32$.
In the figure, we also report the known (asymptotic) lower bounds on the system time (with 3 agents), as derived in~\cite{Bertsimas.vanRyzin:91}, and the system time obtained with the proposed policies in a single-agent scenario.

The performance of both proposed policies is close to optimal for small $\lambda$, as expected. The sensor-based policy behaves well over a large range of target generation rates; in fact, the numerical results suggest that the policy provides a system time that is a constant-factor approximation of the optimum, by a factor of approximately 1.6.

However, as $\lambda$ increases, the performance of the no-communication policy degrades significantly, almost approaching the performance of a single-vehicle system over an intermediate range of values of $\lambda$.  Our intuition in this phenomenon is the following. As agents do not return to their own reference points between visiting successive targets, their efficiency decreases since they are no longer able to effectively separate regions of responsibility. In practice---unless they communicate and concentrate on their Voronoi region, as in the sensor-based policy---agents are likely to duplicate efforts as they pursue the same target, and effectively behave as a single-vehicle system.
Interestingly, this efficiency loss seems to decrease for large $\lambda$, and the numerical results suggest that the no-communication policy recovers a similar performance as the sensor-based policy in the heavy load limit. Unfortunately, we are not able at this time to provide a rigorous analysis of the proposed policies for general values of the target generation rate.

\section{Conclusions}
\label{conclusions}
In this paper we considered two very simple strategies for multiple vehicle routing in the presence of dynamically-generated targets, and analyzed their performance in light load conditions, i.e., when the target generation rate is very small. The strategies we addressed in this paper are based on minimal assumptions on the ability of the agents to exchange information: in one case they do not explicitly communicate at all, and in the other case, agents are only aware of other agents' current location. A possibly unexpected and striking results of our analysis is the following: the collective performance of the agents using such minimal or no-communication strategies is (locally)  optimal, and is in fact as good as that achieved by the best known decentralized strategies. Moreover, the proposed strategies do not rely on the knowledge of the target generation process, and makes minimal assumptions on the target spatial distribution; in fact, the convexity and boundedness assumptions on the support of the spatial distribution can be relaxed, as long as path connectivity of the support, and absolute continuity of the distribution are ensured. Also, the distribution needs not be constant: Indeed, the algorithm will provide a good approximation to a local optimum for the cost function as long as the characteristic time it takes for the target generation process to vary significantly is much greater than the relaxation time of the algorithm. In summary, the proposed strategies can be seen as a learning mechanism in which the agents learn the target generation process, and the ensuing target spatial distribution, and adapt their own behavior to it.

The proposed strategies are very simple to implement, as they only require storage of the coordinates of points visited in the past and simple algebraic calculations; the ``sensor-based" strategy also require a device to estimate the position of other agents. Simple implementation and the absence of active communication makes the proposed strategies attractive, for example, in embedded systems and stealthy applications. The game-theoretic interpretation of our results also provides some insight into how territorial, globally optimal behavior can arise in a population of selfish but rational individuals even without explicit mechanisms for territory marking and defense.

While we were able to prove that the proposed strategies perform efficiently for small values of the target generation rate, little is known about their performance in other regimes. In particular, we have shown numerical evidence that suggests that the first strategy we introduced, requiring no communication, performs poorly when targets are generated very frequently, whereas the performance of the sensor-based strategy is in fact comparable to that of the best known strategies for the heavy load case.

Extensions of this work will include the analysis and design of efficient strategies for general values of the target generation rate, for different vehicle dynamics models (e.g., including differential constraints on the motion of the agents), and heterogeneous systems in which both service requests and agents can belong to several different classes with different characteristics and abilities.

{\bf Acknowledgement}
The research in this paper was inspired by discussion with Dr. Mahbub Gani,  and was performed while the authors were with the Department of Mechanical and Aerospace Engineering at the University of California, Los Angeles.

\bibliographystyle{unsrt}
\bibliography{main,frazzoli}

\end{document}